\RequirePackage{fix-cm}
\documentclass[smallextended]{svjour3}       
\smartqed  
\usepackage{graphicx}
\usepackage{paralist}
\usepackage[misc]{ifsym}
\usepackage{epsfig} 
\usepackage{epstopdf} 
\usepackage{subfigure}
\usepackage[linesnumbered,ruled,vlined]{algorithm2e}
\usepackage[ruled,vlined]{algorithm2e}

\usepackage{makecell}
\usepackage{bm}
\usepackage{cite}
\usepackage{amsmath,amssymb,amsfonts}%
\usepackage{mathrsfs}%
\usepackage[colorlinks=true]{hyperref}

\newtheorem{Example}{Example}

\begin{document}

\title{Numerical method for the inverse scattering problem of acoustic-elastic interaction by random periodic structures
}

\titlerunning{Numerical method for the inverse scattering by random periodic structures}        

\author{Yi Wang         \and
        Lei Lin  \and
        Junliang Lv 
}

\institute{Yi Wang \at
 School of Mathematics, Jilin University, Changchun, 130012, China \\
              \email{wangyi173@163.com}         
            \and
            Lei Lin \at
            School of Mathematics, Jilin University, Changchun, 130012, China\\
             \email{linlei@amss.ac.cn}
             \and
            Junliang Lv \at
            Corresponding author. School of Mathematics, Jilin University, Changchun, 130012, China\\
             \email{lvjl@jlu.edu.cn}
}

\date{Received: date / Accepted: date}

\maketitle

\begin{abstract}
Due to manufacturing defects or wear and tear, industrial components may have uncertainties. In order to evaluate the performance of machined components, it is crucial to quantify the uncertainty of the scattering surface.
This brings up an important class of inverse scattering problems for random interface reconstruction.
In this paper, we present an efficient numerical algorithm for the inverse scattering problem of acoustic-elastic interaction with random periodic interfaces.
The proposed algorithm combines the Monte Carlo technique and the continuation method with respect to the wavenumber, which can accurately reconstruct the key statistics of random periodic interfaces from the measured data of the acoustic scattered field.
In the implementation of our algorithm, a key two-step strategy is employed:
Firstly, the elastic displacement field below the interface is determined by Tikhonov regularization based on the dynamic interface condition;
Secondly, the profile function is iteratively updated and optimised using the Landweber method according to the kinematic
interface condition.
Such a algorithm does not require \emph{a priori} information about the stochastic structures and performs well for both stationary Gaussian and non-Gaussian stochastic processes.
Numerical experiments demonstrate the reliability and effectiveness of our proposed method.

\end{abstract}

\section{Introduction}\label{intro}
Direct and inverse scattering problems form an important research area in modern science and engineering \cite{CK1998book, JxLiLvZh2017ESAIM, JxLiLvZh2017JSC, JxLiLvZh2018CMS, JxLiLvWWZ2022IMA, LinLvLi2024ACM, LinLvNiu2024ACM}.
The scattering theory in periodic structures, also known as gratings in diffractive optics, has great significance in the design and manufacture of optical components, such as improving the sensitivity of optical sensors, increasing the efficiency of beam splitters, and enhancing the performance of corrective lenses and antireflective coatings \cite{BB2001AMO,BLLv2013JOSAA, BL2022book, D1993EJAM, ES1988JCP, NLvG2024IP, ZhLv2018CICP}.
The direct scattering of acoustic-elastic interaction by periodic structures has been extensively investigated in recent years, both in terms of mathematical theories and numerical methods.
A general mathematical
framework for the acoustic-elastic interaction problem in periodic structures was established by variational methods in \cite{HRY2016NMPDE}.
Several computational approaches have also been developed
to numerically solve these problems, such as T matrix method \cite{BS1993JASM} and finite element method \cite{HRY2016NMPDE}.

Due to the dynamic and kinematic interface conditions, the inverse scattering of acoustic-elastic interaction model is much more complex than the inverse single field scattering problem.
There are a few results about theoretical analysis and numerical algorithms of the inverse acoustic-elastic interaction.
The uniqueness results for the incidence of plane waves and point sources can be respectively found in \cite{HKY2016IPIper} and \cite{CuiQW2023JIIP}.
As for the numerical method, the factorization method has been developed for the reconstruction of the bi-periodic interface\cite{HKY2016IPIper}.
So far, the quantitative reconstruction methods for the acoustic-elastic interaction in periodic structures have been little investigated.

The random inverse problems are a class of inverse problems with uncertainty. Since the random inverse problems are closely relevant to the practical application, it is necessary to develop the theoretical analysis and numerical methods for solving such inverse problems. The main challenges of the random inverse problems are twofold: (1) the nonlinearity and ill-posedness of all inverse problems; (2) the complex statistical properties introduced by randomness. Since the complete characterization of the stochastic parameters relies on the second-order moments (e.g., mean, covariance) and higher-order statistics, rather than a single realization of samples, the computational complexity of random inverse problems will significantly increase.
For inverse random source problems, some efficient reconstruction methods for acoustic and elastic waves have been  established in \cite{BCL2016SIAMJUQ, BCLZ2014MC, BCL2017SIAMJNA,LHL2020CIPDE}.
The inverse random potential problem for elastic wave was considered in \cite{LLW2023MMS, LLW2022SIAMJMA}.
As for inverse random interfaces, a class of numerical method has been developed in \cite{BL2021CMR, BLX2020SIAM, WLvL2025IPI} and \cite{GXY2024JCP} for reconstructing electromagnetic and elastic scattering with periodic interfaces, respectively.
However, the research on the inverse scattering problem of acoustic-elastic interaction by random periodic structures has not been investigated yet.

In this paper, we focus on the inverse scattering of acoustic-elastic interaction in random periodic structures. The present work is a nontrivial extension of the method for the inverse electromagnetic scattering problem of random penetrable periodic structures in \cite{WLvL2025IPI}. For such an acoustic-elastic interaction problem, we need to impose the kinematic and kinetic interface conditions on the fluid-solid interface to ensure the continuity of the velocity normal component and the traction, respectively. Due to the coupling of different physical fields and the complex transmission conditions, the original model of the acoustic-elastic interaction is much more complicated than those problems with single wave field. In our method, we firstly reconstruct the surface of each sample accurately, and then derive the key statistical properties of the random structure by Monte Carlo approximation. To reconstruct each sample interface, the elastic displacement field below the interface is first obtained by the Tikhonov regularization method through the dynamic interface condition. Then the profile function is updated by the Landweber iteration method according to the kinematic interface condition.
It is worth pointing out that the grating profiles corresponding to some 'bad' samples can oscillate very violently. If only a single frequency plane wave is used as the incident wave, the reconstruction results cannot capture the real profiles well.
In order to address this issue, we utilize measurements at multiple frequencies and incidence angles.
In addition, to reduce the computational effort, we adopt a modified two-step Monte Carlo continuation algorithm, which uses a small number of samples at low frequencies and can achieve satisfactory computational accuracy.

The rest of the paper is organized as follows. In Section 2, we introduce the mathematical model of the acoustic-elastic interaction in random periodic structure. Section 3 presents the two-step Monte Carlo continuation (TS-MCC) algorithm.
Some convincing numerical experiments are provided in Section 4 to demonstrate the accuracy and reliability of the proposed method. Conclusions and future work are summarized in Section 5.
\section{Mathematical formulations}\label{sec:2}
\begin{figure}[htp]
\begin{center}
  \includegraphics[width=0.5\textwidth]{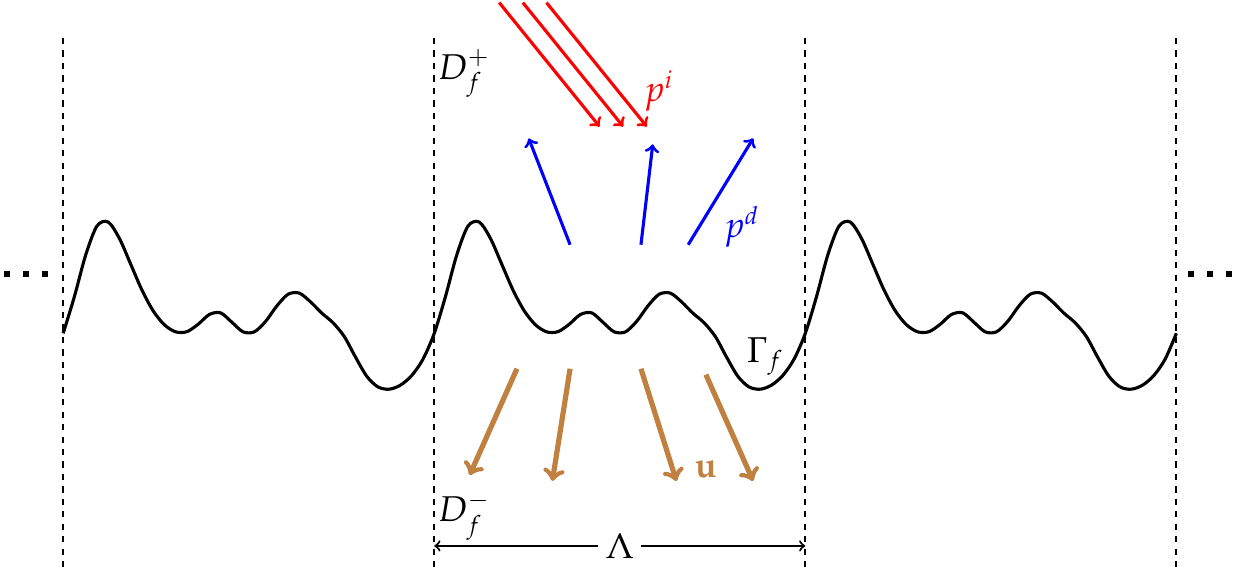}
  \caption{Problem geometry of acoustic-elastic interaction.}\label{problem_geometry}
  \end{center}
\end{figure}

Consider a random acoustic-elastic interaction scattering problem with periodicity in the $x-$axis direction. Assume the random interface is invariant in the $z-$axis direction. Due to the periodicity of the structure, we restrict the system to one periodic unit and define the random periodic interface representing the acoustic-elastic interface as
\begin{align*}
\Gamma_{f}:=\{(x,y)\in {\mathbb{R}^{2}} ; \ x \in [0,\Lambda], y=f(w;x), w \in \Omega \},
\end{align*}
where $w \in \Omega$ denotes a random sample in the probability space $(\Omega, \mathcal{F}, \mu)$, and $f : \Omega \times \mathbb{R} \rightarrow \mathbb{R}$ represents a stationary stochastic process, which can be expressed as the sum of a deterministic function $\tilde{f}(x)$ and a stationary stochastic process $\mathcal{P}(x)$ with zero mean. Moreover, the periods of both the deterministic function $\tilde{f}(x)$ and the stationary stochastic process $\mathcal{P}(x)$ are assumed to be $\Lambda$, i.e., $f(w;x+\Lambda)=f(w;x)$.
As shown in Fig. \ref{problem_geometry}, the domains above and below $\Gamma_{f}$ are respectively defined as
\begin{align*}
D_{f}^{+}:&=\{(x,y)\in {\mathbb{R}^{2}} ; \ x \in [0,\Lambda], y>f(w;x),w \in \Omega \},\\
D_{f}^{-}:&=\{(x,y)\in {\mathbb{R}^{2}}; \ x \in [0,\Lambda], y<f(w;x),w \in \Omega \}.
\end{align*}
Assume that $D_{f}^{+}$ is filled with a homogeneous compressible inviscid fluid with constant mass density $\rho_f > 0$, and $D_{f}^{-}$ is occupied by an isotropic homogeneous elastic solid characterized by the real-valued constant mass density $\rho >0$ and the Lam\'{e} constants $\mu,\lambda \in \mathbb{R}$ satisfying $\mu >0, \lambda + \mu >0$. In order to simplify the notation, define
\begin{align*}
D_{f}:=\{(x,y)\in {\mathbb{R}^{2}} ; \ x \in [0,\Lambda], y\neq f(w;x),w \in \Omega \}.
\end{align*}
\subsection{The acoustic wave}
Suppose that a time-harmonic plane wave illuminates the random interface $\Gamma_f$ from the above,
which induces scattering and transmission phenomena at the acoustic-elastic interface.
Let
\begin{align*}
p^{i}=e^{\mathrm{i}\alpha x -\mathrm{i}\beta y}
\end{align*}
be the incident wave,
where $\alpha=\kappa\sin\theta, \beta=\kappa\cos\theta$, $\theta \in (-\pi / 2, \pi / 2)$ is the angle of incidence. The incident wave satisfies the Helmholtz equation
\begin{align*}
\Delta p^i + \kappa^{2}p^i =0 \quad \mathrm{in} \ D_{f}^{+},
\end{align*}
where $\kappa = \omega / c$ is the wavenumber, $c$ is the speed of sound in the fluid, and $\omega >0$ is the angular frequency. Similarly, for the acoustic total field $p^t$ consisting of the incident field $p^i$ and the diffractive field $p^d$, we have
\begin{align*}
\Delta p^t(w;\cdot) + \kappa^{2}p^t(w;\cdot) =0 \quad \mathrm{in} \ \Omega\times D_{f}^{+}.
\end{align*}
In addition, for the purpose of uniqueness, it is common to seek the quasi-periodic solution which satisfies
\begin{align*}
 p^t(w;x+ \Lambda,y) = e^{\mathrm{i} \alpha \Lambda} p^t(w;x,y)   \quad \mathrm{in} \ \Omega\times D_f^+.
\end{align*}
For each sample $w$, the Rayleigh expansion of the diffractive field $p^d$ can be defined by
\begin{align*}
 p^{d}(w;\cdot)=
 \sum\limits_{n \in \mathbb{Z}} p_{n}(w)e^{\mathrm{i} \alpha_{n} x+ \mathrm{i} \beta_{n}y} \quad  \mathrm{in} \ \Omega\times D_{f}^{+},
\end{align*}
where
\begin{align*}
\alpha_{n}=\alpha + \frac{2n\pi}{\Lambda}, \quad
\beta_{n}=
\begin{cases}
\sqrt{\kappa^{2}-\alpha_{n}^{2}}, &\  |\alpha_{n}| < |\kappa|,
\\
\mathrm{i} \sqrt{\alpha_{n}^{2}-\kappa^{2}}, &\ |\alpha_{n}| > |\kappa|,
\end{cases}
\end{align*}
and $p_{n}(w)$ are Rayleigh coefficients.
To avoid the Wood's anomalies, we assume that $\kappa \neq |\alpha_n|$ for all $ n \in \mathbb{Z}$ \cite{BL2022book}.
\subsection{The elastic wave}

When the acoustic wave hits the surface of the elastic body, an elastic wave is excited inside the elastic body.
The elastic wave satisfies the Navier equation
\begin{align}\label{2.2.1}
\mu \Delta \bm{u}(w;\cdot) + (\lambda +\mu)\nabla \nabla\cdot \bm{u}(w;\cdot) + \omega^2 \rho\bm{u}(w;\cdot) =0 \quad \mathrm{in} \ \Omega\times D_{f}^{-},
\end{align}
where $\bm{u}=(u_1,u_2)^{\top}$ is the displacement field of the elastic wave. Define the traction operator
\begin{align*}
\bm{Tu} := 2\mu \partial_{\bm{n}}\bm{u} + \lambda \bm{n}\nabla\cdot \bm{u}+ \mu
\begin{bmatrix}
  n_2 (\frac{\partial u_2}{\partial x} - \frac{\partial u_1}{\partial y}) \\
  n_1 (\frac{\partial u_1}{\partial y} - \frac{\partial u_2}{\partial x})
\end{bmatrix},
\end{align*}
where $\bm{n}=(n_1, n_2)$ denotes the unit normal vector to the interface $\Gamma_f$ directed into $D_{f}^{+}$.
To ensure the continuity of traction, the following dynamic interface condition is needed
\begin{align*}
-p^t (w;\cdot) \bm{n}= \bm{Tu}(w;\cdot) \quad \mathrm{on} \ \Omega\times \Gamma_f.
\end{align*}
Additionally, the following kinematic interface condition is also required to ensure the continuity of the normal component of the velocity
\begin{align*}
\partial_{\bm{n}} p^t (w;\cdot)= \rho_f \omega^2 \bm{u}(w;\cdot) \cdot \bm{n} \quad \mathrm{on} \ \Omega\times \Gamma_f.
\end{align*}
Similarly, to ensure uniqueness, seek the quasi-periodic solution satisfying
\begin{align*}
\bm{u}(w;x+ \Lambda,y) = e^{\mathrm{i} \alpha \Lambda} \bm{u}(w;x,y)   \quad \mathrm{in} \ \Omega\times D_f^-.
\end{align*}
\subsection{Model problem}

Choose two constants $b^{\pm} \in \mathbb{R}$, such that
\begin{align*}
b^{-} < \min\limits_{w \in \Omega, x \in [0,\Lambda]}f(w;x) \leq \max\limits_{w \in \Omega, x \in [0,\Lambda]}f(w;x) < b^{+}.
\end{align*}
Define
\begin{align*}
\Gamma_{b}^{\pm}&:= \{(x,y)\in {\mathbb{R}^{2}} ; \ x \in [0,\Lambda], y=b^{\pm}\}.
\end{align*}
Under the assumption of small amplitude oscillations in both the solid and fluid,
the direct scattering problem for the random acoustic-elastic interaction periodic structures is formulated as: Find the acoustic field $p^t(w;\cdot)$ and the elastic field $\bm{u}(w;\cdot)$ satisfying
\begin{align*}
\begin{cases}
\Delta p^t(w;\cdot) + \kappa^{2}p^t(w;\cdot) =0 \quad & \mathrm{in} \ \Omega\times D_{f}^{+},
\\
\mu \Delta \bm{u}(w;\cdot) + (\lambda +\mu)\nabla \nabla\cdot \bm{u}(w;\cdot) + \omega^2 \rho\bm{u}(w;\cdot) =0 \quad & \mathrm{in} \ \Omega\times D_{f}^{-},
\\
-p^t (w;\cdot) \bm{n}= \bm{Tu}(w;\cdot) &  \mathrm{on} \ \Omega\times \Gamma_{f},\\
p^t(w;x+ \Lambda,y) = e^{\mathrm{i} \alpha \Lambda} p^t(w;x,y)   \quad &\mathrm{in} \ \Omega\times D_f^+,
\\
\partial_{\bm{n}} p^t (w;\cdot)= \rho_f \omega^2 \bm{u}(w;\cdot) \cdot \bm{n} & \mathrm{on} \ \Omega\times \Gamma_{f},
\\
\bm{u}(w;x+ \Lambda,y) = e^{\mathrm{i} \alpha \Lambda} \bm{u}(w;x,y)   \quad &\mathrm{in} \ \Omega\times D_f^-.
\end{cases}
\end{align*}

This paper mainly focuses on the inverse scattering problem of random acoustic-elastic interaction periodic structures, which can be specifically described as follows: For each sample $\omega$, given the incident field $p^{i}$, determine key statistical quantities of the random interface $y=f(w;x)$ based on the measurements of the diffracted field $p^{d}(\omega;x,b^{+})$ on the straight line $\Gamma_{b}^{+}$.

\section{Reconstruction method}

\subsection{The Helmholtz decomposition}
To begain with, we introduce the Helmholtz decomposition for the displacement field $\bm{u}$. Denote the scalar and vector curl operators
\begin{align*}
\mathrm{curl}  \bm{u} = \partial_{x} u_2 -\partial_{y} u_1, \ \bm{\mathrm{curl}} u= (\partial_{y} u, -\partial_{x} u)^{\top},
\end{align*}
for any vector function $\bm{u}$ and scalar function $u$, respectively.
For each given sample, the Helmholtz decomposition of the displacement field $\bm{u}(w;x,y)$ follows as
\begin{align}\label{3.1.1}
\bm{u}(w;\cdot)= \nabla p_{1}(w;\cdot)+ \bm{\mathrm{curl}}p_{2}(w;\cdot) \quad \mathrm{in} \ \Omega\times D_f^-,
\end{align}
where $p_j(w;\cdot), j=1,2$ are scalar functions.
Substituting \eqref{3.1.1} into \eqref{2.2.1} yields
\begin{align*}
\nabla \big((\lambda + 2\mu)\Delta p_1(w;\cdot)+ \omega^2 \rho p_1(w;\cdot)\big)+\bm{\mathrm{curl}}(\mu \Delta p_2(w;\cdot) +\omega^2 \rho p_2(w;\cdot))=0.
\end{align*}
It is easy to verify $p_j(w;\cdot), j=1,2$ satisfy the Helmholtz equation
\begin{align}\label{3.1.2}
\Delta p_j (w;\cdot)+ {(\kappa_j^-)}^2  p_j(w;\cdot) =0 \quad \mathrm{in} \  \Omega\times D_{f}^{-},
\end{align}
where $\kappa_1^- = \omega  \sqrt{\rho/(\lambda + 2 \mu)}$ denotes the compression  wavenumber, $\kappa_2^- = \omega  \sqrt{\rho/\mu}$  denotes the shear wavenumber.
It follows from the uniqueness that $p_j(w;\cdot),j=1,2$ are quasi-periodic. Thus, we have
\begin{align}\label{3.1.3}
 p_j(w;\cdot)=
 \sum\limits_{n \in \mathbb{Z}} \hat{p}_{jn}(w;y)e^{\mathrm{i} \alpha_{n} x} \quad  \mathrm{in} \ \Omega\times D_{f}^{-},
\end{align}
where
\begin{align*}
 \hat{p}_{jn}(w;y) = \frac{1}{\Lambda}\int_{0}^{\Lambda} p_j(w;x,y)e^{-\mathrm{i} \alpha_{n} x} \mathrm{d}x.
\end{align*}
Substituting \eqref{3.1.3} into \eqref{3.1.2} gives
\begin{align}\label{3.1.4}
\partial_{yy}^{2}\hat{p}_{jn}(w;y) + \beta_{jn}^2\hat{p}_{jn}(w;y)=0 \quad  \mathrm{for} \ y<f(w;x),
\end{align}
where
\begin{align*}
\beta_{jn}=
\begin{cases}
\sqrt{{\kappa_j^-}^{2}-\alpha_{n}^{2}}, &\  |\alpha_{n}| < |\kappa_j^-|,
\\
\mathrm{i} \sqrt{\alpha_{n}^{2}-{\kappa_j^-}^{2}}, &\ |\alpha_{n}| > |\kappa_j^-|.
\end{cases}
\end{align*}
Similarly, we assume $\kappa_j^- \neq |\alpha_n|$ for all $j=1,2$ and $ n \in \mathbb{Z}$.

Using the bounded property of outward propagating waves, it is possible to solve \eqref{3.1.4} analytically and then obtain the explicit solution of \eqref{3.1.2}:
\begin{align*}
p_j(w;x,y) =  \sum_{n \in \mathbb{Z}}\hat{p}_{jn}(w;b^-)e^{\mathrm{i} \alpha_{n} x- \mathrm{i} \beta_{jn}(y-b^-)}  \quad  \mathrm{for} \ y< b^-.
\end{align*}

\subsection{Scattered data}
Assume that the random interface $\Gamma_{f}$ lies between two straight lines
\begin{align*}
\Gamma_{a}^{\pm}:= \{(x,y)\in {\mathbb{R}^{2}}; \ b^{-}<a^{-}<f(\omega;x)<a^{+}<b^{+}, x \in [0,\Lambda], w \in \Omega \}.
\end{align*}
For each given sample $w$, the diffractive field $p^d$ and the scalar potential functions $p_j,j=1, 2$ can be expressed in terms of the single-layer potential, i.e.,
\begin{align*}
p^d(w;x,y)=\int_{0}^{\Lambda} \phi^{+}(w;s)G^{+}(x,y;s,a^{-}) \mathrm{d}s & \quad  \mathrm{in} \  D_{f}^{+}, 
\\
p_j(w;x,y)=\int_{0}^{\Lambda}\phi_j(w;s)G^{-}_j(x,y;s,a^{+}) \mathrm{d}s & \quad  \mathrm{in} \  D_{f}^{-},
\end{align*}
where $\phi^{+},\phi_j \in L^{2}(\Omega;L(0,\Lambda)),j=1,2$ are the unknown quasi-periodic density functions, and $G^{+},G^{-}_j$ are the quasi-periodic Green functions, which are analytically given by
\begin{align}
G^{+}(x,y;s,t)&=\frac{\mathrm{i}}{2 \Lambda}\sum_{n\in\mathbb{Z}}\dfrac{1}{\beta_{n}}e^{\mathrm{i}\alpha_{n}(x-s)+\mathrm{i}\beta_{n}|y-t|}, \quad (x,y)\neq(s,t), \label{3.2.3}\\
G^{-}_j(x,y;s,t)&=\frac{\mathrm{i}}{2 \Lambda}\sum_{n\in\mathbb{Z}}\dfrac{1}{\beta_{jn}}e^{\mathrm{i}\alpha_{n}(x-s)+\mathrm{i}\beta_{jn}|y-t|}, \quad (x,y)\neq(s,t).\label{3.2.4}
\end{align}
Hence, we have
\begin{align}\label{3.2.5}
p^{d}(w;x,b^{+})= \int_{0}^{\Lambda} \phi^{+}(w;s)G^{+}(x,b^{+};s,a^{-})\mathrm{d}s.
\end{align}
Since the scattered field is quasi-periodic, it admits the following expansion
\begin{align}\label{3.2.6}
p^{d}(w;x,b^{+})=\sum_{n\in\mathbb{Z}} p_{n}^{d}(w) e^{\mathrm{i} \alpha_{n} x},
\end{align}
where
\begin{align*}
p_{n}^{d}(w)= \frac{1}{\Lambda}\int_{0}^{\Lambda} p^{d}(w;x,b^{+}) e^{- \mathrm{i} \alpha_{n} x}\mathrm{d}x.
\end{align*}
Similarly, the quasi-periodicity of $\phi^{+}$ and $\phi_j$, $j=1,2$, implies
\begin{align}
\phi^{+}(w ; s)=\sum_{n\in\mathbb{Z}} \phi_{n}^{+}(w) e^{\mathrm{i} \alpha_{n} s},\label{3.2.7}\\
\phi_j(w ; s)=\sum_{n\in\mathbb{Z}} \phi_{jn}(w) e^{\mathrm{i} \alpha_{n} s}.
\label{3.2.8}
\end{align}
It follows from \eqref{3.2.5}-\eqref{3.2.7}  that
\begin{align}\label{3.2.9}
\phi_{n}^{+}(w)=-2\mathrm{i}  \beta_{n}p_{n}^{d}(w)e^{-\mathrm{i} \beta_{n} (b^{+}-a^{-})}.
\end{align}
For simplicity, denote $\partial_1:=\partial_x, \partial_2 := \partial_y, \partial_{11}:=\partial_{xx}, \partial_{12}:=\partial_{xy}, \partial_{21}:=\partial_{yx}, \partial_{22}:=\partial_{yy}$.
Define eight operators $T_{f}^{+}, T_{i}^{+}, T_{i}, T_{ik}, :L^{2}(\Omega;L^{2}(0,\Lambda))\rightarrow L^{2}(\Omega;L^{2}(0,\Lambda)),$ $(i,k=1,2)$ by
\begin{align}
(T_{f}^{+} \phi^{+})(w;x)&=\int_{0}^{\Lambda}\phi^{+}(w;s)G^{+}(x,f(w;x);s,a^{-}) \mathrm{d}s, \label{3.2.10}\\
(T_{i}^{+} \phi^{+})(w;x)&=\int_{0}^{\Lambda}\phi^{+}(w;s)\partial_i G^{+}(x,f(w;x);s,a^{-}) \mathrm{d}s, \label{3.2.11}\\
(T_{i} \phi_j)(w;x)&=\int_{0}^{\Lambda}\phi_j(w;s)\partial_i G^{-}_j(x,f(w;x);s,a^{+}) \mathrm{d}s, \label{3.2.12}\\
(T_{ik} \phi_j)(w;x)&=\int_{0}^{\Lambda}\phi_j(w;s)\partial_{ik} G^{-}_j(x,f(w;x);s,a^{+}) \mathrm{d}s. \label{3.2.13}
\end{align}
Substituting \eqref{3.2.3}, \eqref{3.2.4} and  \eqref{3.2.7}-\eqref{3.2.9} into \eqref{3.2.10}-\eqref{3.2.13}   yields
\begin{align}
(T_{f}^{+} \phi^{+})(w;x)&=\sum_{n\in\mathbb{Z}} \psi_{n}^{+}(w) e^{\mathrm{i} \alpha_{n} x + \mathrm{i} \beta_{n} f(w;x)}, \label{3.2.14}\\
\bm{K}^+ \phi^{+}:= (T_{1}^{+} \phi^{+}, T_{2}^{+} \phi^{+} )(w;x)&=\sum_{n\in\mathbb{Z}} (\mathrm{i} \alpha_n,\mathrm{i}\beta_n)\psi_{n}^{+}(w) e^{\mathrm{i} \alpha_{n} x + \mathrm{i} \beta_{n} f(w;x)},\label{3.2.15}\\
\bm{K}^-\phi_j:=(T_{1} \phi_{j}, T_{2} \phi_{j})(w;x)&=\sum_{n\in\mathbb{Z}} (\mathrm{i} \alpha_n,-\mathrm{i}\beta_{jn})\psi_{jn}(w) e^{\mathrm{i} \alpha_{n} x - \mathrm{i} \beta_{jn} f(w;x)},\label{3.2.16}\\
\bm{H}^-\phi_j:=
\begin{bmatrix}
  T_{11} \ \   T_{12}\\
  T_{21} \ \  T_{22}
\end{bmatrix}
\phi_j(w;x)&=\sum_{n\in\mathbb{Z}}
\begin{bmatrix}
  -\alpha_n^2 \ \ \ \    \alpha_n \beta_{jn}\\
  \alpha_n \beta_{jn} \ \  -\beta_{jn}^2
\end{bmatrix}
\psi_{jn}(w) e^{\mathrm{i} \alpha_{n} x - \mathrm{i} \beta_{jn} f(w;x)},\label{3.2.17}
\end{align}
where $\psi_{n}^{+}(w)= p_{n}^{d}(w)e^{-\mathrm{i} \beta_{n} b^{+}}, \psi_{jn}(w)= \frac{\mathrm{i}}{2 \beta_{jn}}\phi_{jn}(w)e^{\mathrm{i} \beta_{jn} a^{+}}$.
In practice, we use the regularization (see also \cite{EHN1996book}) to restrict the exponential growth of noise for $b^+>0$, i.e.,
\begin{align*}
\psi_{n}^{+}(\omega)=
\begin{cases}
p_{n}^{d}(\omega)e^{-\mathrm{i} \beta_{n} b^{+}}, &\  \mathrm{for} \  \kappa^{+}>|\alpha_{n}|,
\\
p_{n}^{d}(\omega)\dfrac{e^{\mathrm{i} \beta_{n} b^{+}}}{e^{2\mathrm{i} \beta_{n} b^{+}}+\gamma}, &\  \mathrm{for} \  \kappa^{+}<|\alpha_{n}|,
\end{cases}
\end{align*}
where $\gamma$ is some positive regularization parameter.

\subsection{The nonlinear problem}

For any $w \in \Omega$, according to the dynamic interface condition, we get
\begin{align}\label{3.3.1}
\begin{split}
-\big((T_{f}^{+} \phi^{+})(w;x)+p^i(w;x)\big)\bm{n} =
& 2 \mu \big(\bm{H}^{-} \phi_1 +\bm{P} \bm{H}^{-} \phi_2 \big )
 \bm{n}\\[2mm]
& +\lambda  \text{tr}(\bm{H}^{-}) \phi_1 \bm{n}
+\mu \text{tr}(\bm{H}^{-} )\phi_2 \widetilde{\bm{I}}\bm{ P^{\top}}\bm{n},
\end{split}
\end{align}
where $\widetilde{\bm{I}}$ is a unitary matrix of second order and $\bm{P}=\begin{bmatrix}
  \ 0 \ \ \    1\\
  -1 \ \ \  0
\end{bmatrix}$.
Substituting \eqref{3.2.14}, \eqref{3.2.17} into  \eqref{3.3.1} yields
\begin{align*}
&-\Big(\sum_{n\in\mathbb{Z}} \psi_{n}^{+}(w) e^{\mathrm{i} \alpha_{n} x + \mathrm{i} \beta_{n} f(w;x)}
+ e^{\mathrm{i} \alpha x-\mathrm{i} \beta f(w;x)}\Big)\bm{n} \\
=& 2\mu \sum_{n\in\mathbb{Z}}
\begin{bmatrix}
  (-\alpha_n^2, \alpha_n \beta_{2n})\bm{\chi}_n & (\alpha_n \beta_{1n},\ -\beta_{2n}^2)\bm{\chi}_n\\
  (\alpha_n \beta_{1n},\ \alpha_n^2)\bm{\chi}_n & -(\beta_{1n}^2,\ \  \alpha_n \beta_{2n})\bm{\chi}_n
\end{bmatrix}
 \bm{n} \\
& + \lambda  \sum_{n\in\mathbb{Z}}  \big(-\alpha_n^2-\beta_{1n}^2, \ 0\big)\bm{\chi}_n \bm{n}  +
\mu \sum_{n\in\mathbb{Z}}
\begin{bmatrix}
  0 & (\alpha_n^2+ \beta_{2n}^2)\chi_{2n} \\
 -(\alpha_n^2+ \beta_{2n}^2)\chi_{2n}&  0
\end{bmatrix}\bm{n},
\end{align*}
where
\begin{align*}
\bm{\chi}_n =\begin{bmatrix}
 \chi_{1n} \\
  \chi_{2n}
\end{bmatrix} =\begin{bmatrix}
 \psi_{1n} e^{\mathrm{i} \alpha_{n} x - \mathrm{i} \beta_{1n} f(w;x)}\\
  \psi_{2n} e^{\mathrm{i} \alpha_{n} x - \mathrm{i} \beta_{2n} f(w;x)}
\end{bmatrix}.
\end{align*}
Since the series in the above equation decreases exponentially with respect to $|n|$, we truncate them into a finite number of terms by choosing a sufficiently large $N$, which leads to the approximation equation
\begin{align}\label{3.3.2}
\begin{split}
&-\Big(\sum_{n =-N}^{N} \psi_{n}^{+}(w) e^{\mathrm{i} \alpha_{n} x + \mathrm{i} \beta_{n} f(w;x)}
+ e^{\mathrm{i} \alpha x-\mathrm{i} \beta f(w;x)}\Big)\bm{n} \\
=& 2\mu \sum_{n =-N}^{N}
\begin{bmatrix}
  (-\alpha_n^2,\ \ \alpha_n \beta_{2n})\bm{\chi}_n\ \   (\alpha_n \beta_{1n},\ -\beta_{2n}^2)\bm{\chi}_n\\
  (\alpha_n \beta_{1n},\ \alpha_n^2)\bm{\chi}_n\ \ \ \  -(\beta_{1n}^2,\ \  \alpha_n \beta_{2n})\bm{\chi}_n
\end{bmatrix}
 \bm{n} \\
& + \lambda  \sum_{n =-N}^{N} \big(-\alpha_n^2-\beta_{1n}^2, \ 0\big)\bm{\chi}_n \bm{n}  +
\mu \sum_{n =-N}^{N}
\begin{bmatrix}
  0 & (\alpha_n^2+ \beta_{2n}^2)\chi_{2n} \\
 -(\alpha_n^2+ \beta_{2n}^2)\chi_{2n}&  0
\end{bmatrix}\bm{n}.
\end{split}
\end{align}

In the meanwhile, based on the kinematic interface condition, we have
\begin{align*}
  \Big((\mathrm{i} \alpha, -\mathrm{i} \beta) e^{\mathrm{i} \alpha x-\mathrm{i} \beta f(w;x)} + \bm{K}^+ \phi^{+} \Big) \bm{n} - \rho_f \omega^2 \Big(
T_1 \phi_1 +T_2 \phi_2,
  T_2 \phi_1 -T_1 \phi_2
\Big)
 \bm{n}
 =0.
\end{align*}
Combining the above equation together with \eqref{3.2.15} and \eqref{3.2.16}, we obtain
\begin{align*}
 0=& \bigg(  \alpha e^{\mathrm{i} \alpha x-\mathrm{i} \beta f(w;x)}
 \ +  \sum_{n\in\mathbb{Z}}  \alpha_n \psi_{n}^{+}(w) e^{\mathrm{i} \alpha_{n} x + \mathrm{i} \beta_{n} f(w;x)}\bigg) n_1\\
& + \bigg(- \beta e^{\mathrm{i} \alpha x -\mathrm{i} \beta f(w;x) }
+\sum_{n\in\mathbb{Z}} \beta_{n}\psi_{n}^{+}(w) e^{\mathrm{i} \alpha_{n} x + \mathrm{i} \beta_{n} f(w;x)} \bigg) n_2\\
& -\rho_f \omega^2 \sum_{n\in\mathbb{Z}}\Big( \alpha_{n} \chi_{1n}- \beta_{2n} \chi_{2n}\Big) n_1
 +\rho_f \omega^2  \sum_{n\in\mathbb{Z}}\Big( \beta_{1n} \chi_{1n}+ \alpha_{n} \chi_{2n}\Big) n_2,
\end{align*}
where \begin{align*}
\bm{n}=\begin{bmatrix}
 n_1 \\
  n_2
\end{bmatrix} = \frac{1}{\sqrt{f'(w;x)^{2}+1}}\begin{bmatrix}
 -f'(w;x)\\
1
\end{bmatrix}.
\end{align*}
Here, $f'(w;x)$ is the derivative of $f(w;x)$ with respect to $x$.

Similar to \eqref{3.3.2}, the approximate equation can be truncated into a finite summation
\begin{align}\label{3.3.3}
\begin{split}
 0=& \bigg(  \alpha e^{\mathrm{i} \alpha x-\mathrm{i} \beta f(w;x)}
 \ +  \sum_{n=-N}^{N}  \alpha_n \psi_{n}^{+}(w) e^{\mathrm{i} \alpha_{n} x + \mathrm{i} \beta_{n} f(w;x)}\bigg) n_1\\
& +\bigg(- \beta e^{\mathrm{i} \alpha x -\mathrm{i} \beta f(w;x) }
+\sum_{n=-N}^{N} \beta_{n}\psi_{n}^{+}(w)
e^{\mathrm{i} \alpha_{n} x + \mathrm{i} \beta_{n} f(w;x)} \bigg)n_2\\
& -\rho_f \omega^2 \sum_{n=-N}^{N}\Big( \alpha_{n} \chi_{1n}- \beta_{2n} \chi_{2n}\Big)n_1
+\rho_f \omega^2  \sum_{n=-N}^{N}\Big( \beta_{1n} \chi_{1n}+ \alpha_{n} \chi_{2n}\Big)n_2.
\end{split}
\end{align}
Clearly, \eqref{3.3.2} and \eqref{3.3.3} form a system of nonlinear equations with respect to $f$.

\subsection{The two-step Monte Carlo continuation algorithm}
In order to reconstruct the key statistics of the random structures from measurements of the scattered field on a straight line above the periodic structure, we propose a two-step Monte Carlo continuation (TS-MCC) algorithm.
We employ the Monte Carlo technique to sample the probability space, and let $M \in \mathbb{N^{+}}$ be the sample size. For each sample $w_{m}, m=1,2,\ldots,M$, we denote the reconstructed random periodic interface by $f_{m}(x),m=1,2,\ldots,M$. Without loss of generality, we take the period $\Lambda$ to be $2 \pi$ from now on. Thus, the Fourier series expansion of the sample $f_{m}, m=1,2,\ldots,M$, can be written as
\begin{align*}
f_{m}(x)=a_{m,0}+\sum_{p=1}^{\infty}\left[a_{m,2p-1}\cos(px)+a_{m,2p}\sin(px)   \right].
\end{align*}

It is obvious that the finite term expansion of the above equation is a reasonable approximation of the random periodic interface.
As discussed in \cite{BLX2020SIAM}, the low-frequency data can only recover the approximate rough outline of the target, while the high-frequency data can capture the local geometric detail of the target. Therefore the reconstruction can be significantly improved by using multi-frequency data.
Based on the relationship between frequency and wave number in the acoustic-elastic interaction model, such a multi-frequency inversion algorithm can be implemented by employing a set of appropriately selected wavenumbers.
To describe the steps of the algorithm in detail, we first introduce two parameter vectors
\begin{align}\label{ZZ}
\boldsymbol{\kappa}= [\kappa_{1}, \kappa_{2}, \ldots, \kappa_{Q} ]^{\top}, \ \
\mathbf{Z} = [z_{1}, z_{2}, \ldots, z_{Q} ]^{\top},
\end{align}
where the positive real number $\kappa_{j}, j=1,2,\ldots,Q$, are the wavenumber above the interface and are ordered in an increasing sequence, the non-negative integer $z_{j}, j=1,2,\ldots,Q$, are the largest integer which are smaller than or equal to $\kappa_{j}$. For each sample $w_{m},m=1,2,\ldots,M$, let
\begin{align*}
\mathbf{a}_{m,z_{Q}}= [a_{m,0}, a_{m,1}, \ldots, a_{m,2 z_{Q}} ]^{\top}
\end{align*}
be a vector consisting of the unknowns, which are the Fourier coefficients of the profile function. Initially, we set $a_{m,0}=b^{+}$ and $a_{m,p}=0,p=1,2,\ldots,2 z_{Q}$.

To determine each coefficient $a_{m,p}=0,p=0,1,\ldots,2 z_{Q}$, a two-step approach is applied to the reconstruction process for each $\kappa_{j},j=1,2,\ldots,Q$.
Then, we describe the specific process using $\kappa_{1}$ as an example.


{\bf Step 1}. Approximately compute the elastic displacement field in the solid below the interface. At this point we've got  $\mathbf{a}_{m,z_{1}}=[a_{m,0}, a_{m,1}, \ldots, a_{m,2z_{1}} ]^{\top}$, where $a_{m,0}=b^{+}, a_{m,p}=0,p=1,2,\ldots,2z_{1}$. We use the truncated Fourier series
\begin{align*}
f_{m,z_{1}}(x)=a_{m,0}+\sum_{p=1}^{z_{1}}\left[a_{m,2p-1}\cos(px)+a_{m,2p}\sin(px)   \right]
\end{align*}
to approximate the profile function $f(w_{m};x)$.
In order to obtain a more accurate approximation, we use illuminations from multiple angles $\theta_{l} \in (-\pi/2,\pi/2), l=1,2,\ldots,L$. By denoting $\alpha_{l} = \kappa_{1} \sin(\theta_{l}),\beta_{l}= \kappa_{1}\cos(\theta_{l})$,  and $x_{j} = \frac{\pi}{N'}j, \bm{n}_{j}=(n_{1,j}, n_{2,j})^{\top} = \left(f'_{m,z_{1}}(x_{j}),1\right)^{\top} / \left((f'_{m,z_{1}}(x_{j}))^{2}+1\right)^{1/2}, j=0,1,\ldots,2N'$, we can define
\begin{align*}
A_{11}^{(l)} & =\bigg[ \big(
-{\alpha_n^{(l)}}^2 n_{1,j}+ \alpha_n^{(l)}\beta_{1n}^{(l)}n_{2,j} \big)
e^{\mathrm{i} \alpha_{n}^{(l)}  x_j - \mathrm{i} \beta_{1n}^{(l)}  f_{m,z_{1}}(x_j)}
\bigg]_{(2N'+1)\times (2N+1)},\\
A_{12}^{(l)} & =\bigg[\big(
\alpha_n^{(l)}\beta_{2n}^{(l)} n_{1,j}- {\beta_{2n}^{(l)}}^2 n_{2,j}\big)
e^{\mathrm{i} \alpha_{n}^{(l)}  x_j - \mathrm{i} \beta_{2n}^{(l)}  f_{m,z_{1}}(x_j)}
\bigg]_{(2N'+1)\times (2N+1)},\\
A_{21}^{(l)} & = \bigg[\big(
\alpha_n^{(l)}\beta_{1n}^{(l)} n_{1,j}- {\beta_{1n}^{(l)}}^2 n_{2,j}\big)
e^{\mathrm{i} \alpha_{n}^{(l)}  x_j - \mathrm{i} \beta_{1n}^{(l)}  f_{m,z_{1}}(x_j)}
\bigg]_{(2N'+1)\times (2N+1)},\\
A_{22}^{(l)} & =\bigg[\big(
{\alpha_n^{(l)}}^2 n_{1,j}- \alpha_n^{(l)}\beta_{2n}^{(l)}n_{2,j}\big)
e^{\mathrm{i} \alpha_{n}^{(l)}  x_j - \mathrm{i} \beta_{2n}^{(l)}  f_{m,z_{1}}(x_j)}
\bigg]_{(2N'+1)\times (2N+1)},\\
B_{11}^{(l)} & =\bigg[
\big(
-{\alpha_n^{(l)}}^2 - {\beta_{1n}^{(l)}}^2 \big)n_{1,j}
e^{\mathrm{i} \alpha_{n}^{(l)}  x_j - \mathrm{i} \beta_{1n}^{(l)}  f_{m,z_{1}}(x_j)}
\bigg]_{(2N'+1)\times (2N+1)},\\
B_{21}^{(l)} & =\bigg[
\big(
-{\alpha_n^{(l)}}^2 - {\beta_{1n}^{(l)}}^2 \big)
n_{2,j}
e^{\mathrm{i} \alpha_{n}^{(l)}  x_j - \mathrm{i} \beta_{1n}^{(l)}  f_{m,z_{1}}(x_j)}
\bigg]_{(2N'+1)\times (2N+1)},\\
B_{12}^{(l)} &= B_{22}^{(l)} = C_{11}^{(l)}  = C_{21}^{(l)}=\bigg[ \
0 \
\bigg]_{(2N'+1)\times (2N+1)}, \\
C_{12}^{(l)} & =\bigg[
\big(
{\alpha_n^{(l)}}^2 + {\beta_{2n}^{(l)}}^2 \big) n_{2,j}
e^{\mathrm{i} \alpha_{n}^{(l)}  x_j - \mathrm{i} \beta_{2n}^{(l)}  f_{m,z_{1}}(x_j)}
\bigg]_{(2N'+1)\times (2N+1)},\\
C_{22}^{(l)} & =\bigg[
- \big(
{\alpha_n^{(l)}}^2 + {\beta_{2n}^{(l)}}^2 \big) n_{1,j}
e^{\mathrm{i} \alpha_{n}^{(l)}  x_j - \mathrm{i} \beta_{2n}^{(l)}  f_{m,z_{1}}(x_j)}
\bigg]_{(2N'+1)\times (2N+1)},\\
G_{1}^{(l)} & =\bigg[ -\bigg(e^{\mathrm{i} \alpha^{(l)} x_{j}-\mathrm{i}\beta^{(l)} f_{m,z_{1}}(x_{j}) } + \sum_{n=-N}^{N}{\psi_{n}^{+}}^{(l)} e^{\mathrm{i} \alpha_{n}^{(l)} x_{j}+\mathrm{i} \beta_{n}^{(l)}f_{m,z_{1}}(x_{j})} \bigg)
n_{1,j}
\bigg]_{(2N'+1)\times 1},  \\
G_{2}^{(l)} & =\bigg[ -\bigg(e^{\mathrm{i} \alpha^{(l)} x_{j}-\mathrm{i}\beta^{(l)} f_{m,z_{1}}(x_{j}) } + \sum_{n=-N}^{N}{\psi_{n}^{+}}^{(l)} e^{\mathrm{i} \alpha_{n}^{(l)} x_{j}+\mathrm{i} \beta_{n}^{(l)}f_{m,z_{1}}(x_{j})} \bigg)
n_{2,j}
\bigg]_{(2N'+1)\times 1}.
\end{align*}
It follows from \eqref{3.3.2} that
\begin{align}\label{3.17}
A\bm{\Psi}=G,
\end{align}
where
\begin{align*}
& A=\big[A^{(1)}; \ldots;A^{(L)}\big],
\bm{\Psi}=\big[\bm{\Psi}^{(1)}; \ldots; \bm{\Psi}^{(L)} \big],
G=\big[G^{(1)}; \ldots; G^{(L)}\big],
\\[3mm]
& A^{(l)}=2\mu
\begin{bmatrix}
  A_{11}^{(l)}\ \    A_{12}^{(l)}\\
   A_{21}^{(l)} \ \   A_{22}^{(l)}
\end{bmatrix}
+ \lambda
\begin{bmatrix}
  B_{11}^{(l)}\ \    B_{12}^{(l)}\\
   B_{21}^{(l)} \ \   B_{22}^{(l)}
\end{bmatrix}
+ \mu
\begin{bmatrix}
  C_{11}^{(l)}\ \    C_{12}^{(l)}\\
   C_{21}^{(l)} \ \   C_{22}^{(l)}
\end{bmatrix}, \\[3mm]
& \bm{\Psi}^{(l)}=\begin{bmatrix}
  \bm{\psi}_1^{(l)}\\
   \bm{\psi}_2^{(l)}
\end{bmatrix} =  \begin{bmatrix}
  \bigg[\psi_{1,-N}^{(l)}, \psi_{1,-(N-1)}^{(l)}, \ldots,\psi_{1,(N-1)}^{(l)}, \psi_{1,N}^{(l)} \bigg]^{\top} \\
   \bigg[\psi_{2,-N}^{(l)}, \psi_{2,-(N-1)}^{(l)}, \ldots,\psi_{2,(N-1)}^{(l)}, \psi_{2,N}^{(l)} \bigg]^{\top}
\end{bmatrix}
, \quad
G^{(l)}=
\begin{bmatrix}
  G_{1}^{(l)}\\
   G_{2}^{(l)}
\end{bmatrix}.
\end{align*}

Owing to the ill-posedness of \eqref{3.17}, regularization becomes imperative prior to numerical resolution. A practical methodology involves introducing a stabilization scheme through the selection of an appropriate regularization parameter $\varepsilon$, followed by solving the modified system
\begin{align*}
\big(A+\varepsilon I\big)\bm{\Psi}=G.
\end{align*}
Having obtained the numerical solution $\bm{\Psi}$, the elastic displacement field below the interface can be approximated by
\begin{align*}
\bm{u}(w;x,f(w;x))
& =\sum_{n=-N}^{N}
\begin{bmatrix}
 \mathrm{i} \alpha_n  & - \mathrm{i} \beta_{2n} \\
  -\mathrm{i} \beta_{1n} &  - \mathrm{i} \alpha_{n}
\end{bmatrix}
\begin{bmatrix}
   \chi_{1n} \\
 \chi_{2n}
\end{bmatrix} \\
& =
\sum_{n=-N}^{N}
\begin{bmatrix}
 \mathrm{i} \alpha_n & - \mathrm{i} \beta_{2n} \\
  -\mathrm{i} \beta_{1n} &  - \mathrm{i} \alpha_{n}
\end{bmatrix}
\begin{bmatrix}
   \psi_{1n}e^{\mathrm{i} \alpha_{n} x - \mathrm{i} \beta_{1n} f(w;x)} \\
 \psi_{2n}e^{\mathrm{i} \alpha_{n} x - \mathrm{i} \beta_{2n} f(w;x)}
\end{bmatrix}.
\end{align*}

{\bf Step 2}. Update the reconstructed profile function. Define
\begin{align*}
 &J^{(l)}(\mathbf{a}_{m,z_{1}})   \\
= &\Bigg\|
  \sum_{n =-N}^{N}
\bigg( \alpha_n^{(l)} n_1 + \beta_n^{(l)} n_2\bigg)
 {\psi_{n}^{+}}^{(l)} (w) e^{\mathrm{i} \alpha_{n}^{(l)}  x + \mathrm{i} \beta_{n}^{(l)}  f_{m,z_{1}}(x)}
  + \bigg( \alpha^{(l)}n_1 - \beta^{(l)}n_2\bigg)
e^{\mathrm{i} \alpha^{(l)} x-\mathrm{i} \beta^{(l)} f_{m,z_{1}}(x)} \\
&  \quad- \rho_f \omega_{q_1}^2
\sum_{n=-N}^{N}
\bigg( \alpha_{n}^{(l)} n_1 - \beta_{1n}^{(l)}n_2
\bigg)\chi_{1n}^{(l)}
+\rho_f \omega_{q_1}^2
\sum_{n=-N}^{N}
\bigg( \beta_{2n}^{(l)} n_1 +  \alpha_{n}^{(l)} n_2
\bigg)\chi_{2n}^{(l)}
\Bigg\|_{L^{2}(\Omega;L^{2}(0,\Lambda))}^{2},
\end{align*}
and denote  $\mathbf{J}(\mathbf{a}_{m,z_{1}})=[J^{(1)}(\mathbf{a}_{m,z_{1}}),\ldots,$ $
J^{(L)}(\mathbf{a}_{m,z_{1}})]^{\top}$. Then the system of nonlinear equations \eqref{3.3.3} can be reformulated as
\begin{align}\label{3.19}
\mathbf{J}(\mathbf{a}_{m,z_{1}})=\mathbf{0},
\end{align}
where $\mathbf{J}(\mathbf{a}_{m,z_{1}}): \mathbb{R}^{2 z_{1} +1} \rightarrow \mathbb{R}^{L}$.
To reduce the computational cost and instability, the nonlinear Landweber iteration (see also \cite{EHN1996book})
\begin{align}\label{3.20}
\mathbf{a}_{m,z_{1}}^{(t+1)}=\mathbf{a}_{m,z_{1}}^{(t)}-\eta_{\kappa_{1}}\mathbf{D}\mathbf{J}^{\top}(\mathbf{a}_{m,z_{1}}^{(t)})\mathbf{J}(\mathbf{a}_{m,z_{1}}^{(t)}), \quad t=0,1,2, \ldots,
\end{align}
is used to solve \eqref{3.19}, where $\mathbf{a}_{m,z_{1}}^{(0)}=\mathbf{a}_{m,z_{1}}$, the relaxation parameter $\eta_{\kappa_{1}}$  depends on the wavenumber $\kappa_{1}$, and the Jacobi matrix
\begin{align*}
\mathbf{D}\mathbf{J}=\left(\dfrac{\partial J^{(l)}}{\partial a_{m,j}} \right)_{l=1,2,\ldots,L, j=0,1,\ldots,2z_{1}}.
\end{align*}
After the proper number of iterations, we denote the updated $\mathbf{a}_{m,z_{1}}$ by $\tilde{\mathbf{a}}_{m,z_{1}}=[\tilde{a}_{m,0}, \tilde{a}_{m,1}, \ldots, \tilde{a}_{m,2z_{1}} ]^{\top}$.
Namely, the first $2 z_{1} +1$ terms of $\tilde{a}_{m,2z_{Q}} $ are obtained using the above nonlinear Lundweber method at wavenumber $\kappa_1$.
This concludes the calculation process in Step 2.

The reconstruction result of the first round is usually not accurate enough, as it contains only a few Fourier modes. Next, we increase the wavenumber to $\kappa_{2}$ and start the next round of reconstruction.
Denote $\mathbf{a}_{m,q_{2}}=[a_{m,0}, a_{m,1}, \ldots, a_{m,2z_{2}} ]^{\top}$.
The objective remains analogous: Determine a truncated Fourier expansion
\begin{align*}
f_{m,z_{2}}(x)=a_{m,0}+\sum_{p=1}^{z_{2}}\left[a_{m,2p-1}\cos(px)+a_{m,2p}\sin(px)   \right]
\end{align*}
to approximate the profile function $f(w_{m};x)$. Initially, we set
\begin{align*}
a_{m,p}: =
\begin{cases}
\tilde{a}_{m,p}, &\  0\leq p\leq2z_{1} ,\\
\\
0, &\  2z_{1}< q \leq2z_{2}.\\
\end{cases}
\end{align*}
By repeating Steps 1-2, we can obtain the approximate coefficient vector $\tilde{\mathbf{a}}_{m,z_{2}} =[\tilde{a}_{m,0}, \tilde{a}_{m,1}, \ldots , \tilde{a}_{m,2z_{2}} ]^{\top}$ of the truncated Fourier function corresponding to the wavenumber $\kappa_{2}$.

Repeat the above process until the biggest wavenumber $\kappa_{Q}$ is reached.
Eventually, we get a final approximation of the truncated Fourier coefficients $\tilde{\mathbf{a}}_{m,z_{ Q}} =[\tilde{a}_{m,0}, \tilde{a}_{m,1}, \ldots , \tilde{a}_{m,2z_{Q}} ]^{\top}$.

\subsection{Modified two-step Monte Carlo continuation algorithm}
We now recall the main idea of the two-step Monte Carlo continuation algorithm.
For each sample $w_{m}, m=1,2,\ldots,M$, the initial value of the high-frequency reconstruction is given by the previous low-frequency reconstruction results.
The low-frequency data can only recover the rough outline of the interface, while the high-frequency data can capture the local details of the interface.
Therefore, in order to improve the computational efficiency, we can use fewer samples at low frequencies and choose more samples at high frequencies, i.e., choose different numbers of samples according to different frequencies.
We can modify the two-step Monte Carlo continuation algorithm according to this characteristic, and describe its detailed procedure as follows.

Let $\Omega_M = \{w_m \in \Omega; \ m=1,2,\ldots,M\}$ be the set of samples. For the wavenumber $\kappa_j$, we choose $\Omega_{M_j}=\{w_m \in \Omega_M; \ m=1,2,\ldots,M_j\}$ to reconstruct the grating profile, where $\Omega_{M_j}$ consists of the first $M_j$ samples. For each wavenumber $\kappa_j,j=1,2,\ldots, Q$, we only compute the samples in $\Omega_{M_j}$. In practice, we can take the initial value $f^0=b^{+}$ and $f^{j-1}(x)=\frac{1}{M_{j-1}} \sum_{m=1}^{M_{j-1}}f_{m,z_{j-1}}(x), j=2,3,\ldots, Q$, and reconstruct the profile functions $f_{m ,z_{j}}(x), $ $ m=1,2,\ldots,M_{j}$, using the two-step algorithm corresponding to the wavenumber $\kappa_j$.
\subsection{Computation of statistics}
Recall that the random interface is the sum of the deterministic function $\tilde{f}(x)$ and the zero-mean static random process $\mathcal{P}(x)$.
We can use the mean function
\begin{align*}
\bar{f}(x)=\bar{a}_{0}+\sum_{p=1}^{z_{Q}}\left[\bar{a}_{2p-1}\cos(px)+\bar{a}_{2p}\sin(px) \right]
\end{align*}
as an approximation to $\tilde{f}(x)$, where
\begin{align*}
\bar{a}_{p} = \dfrac{1}{M} \sum_{m=1}^{M} \tilde{a}_{m,p}, \quad p=0,1,\ldots,2z_{Q}.
\end{align*}
In addition, the covariance function can be recovered by the Monte Carlo technique.
As is well known, the covariance function of a stochastic process can be expressed as
\begin{align}\label{3.21}
c(s,t) = \mathbb{E}[f(s)-\mathbb{E}f(s)][f(t)-\mathbb{E}f(t)], \quad s,t \in (0,2\pi).
\end{align}
Then, to obtain the covariance matrix $\bar{C}=(\bar{c}_{ij})$, we choose a finite number of discrete points $\{s_i\}_{i=1}^{n}$ and $\{t_j\}_{j=1}^{n}$ to discretise \eqref{3.21}, where
\begin{align*}
\bar{c}_{ij} = \frac{1}{M}\sum_{m=1}^{M}[f_{m}(s_i)-\bar{f}(s_i)][f_{m}(t_j)-\bar{f}(t_j)], \quad i,j=1,2,\ldots,n.
\end{align*}
The probability distribution of a stationary Gaussian stochastic process is completely determined by its mean and covariance functions.
However, for a non-Gaussian process with non-zero higher-order moments (the third-order skewness and the fourth-order kurtosis, etc.), only the first two moments cannot fully characterize its statistics.
To fully characterize a non-Gaussian process, knowledge of all its higher-order moments or the complete probability density function is theoretically required.
In this paper, we choose the Kernel Density Estimation method \cite{W2018ITMWC} to estimate the probability density function
\begin{align*}
\widehat{f}_{d}(x)=\frac{1}{Md}\sum_{i=1}^{M}K(\frac{x-x_{i}}{d})
\end{align*}
at any location to reconstruct a non-Gaussian process.
Here, the kernel function $K(\cdot)$ is non-negative and satisfies $\int_{-\infty}^{\infty}K(x) \mathrm{d} x =1$, and the bandwidth parameter $d > 0$ denotes the smoothness of the density estimate.
The Kernel Density Estimation is a non-parametric estimation method that does not require \emph{a priori} knowledge.
It constructs density estimates by locally weighted averaging some data points and is capable to handle complex data.
Therefore, it is well suited to deal with our problem.

In conclusion, we outline the complete pseudocode of the modified two-step Monte Carlo continuation algorithm in Algorithm 1.
\begin{algorithm}
    \SetAlgoLined 
	\caption{Modified two-step Monte Carlo continuation algorithm}
   Take $\boldsymbol{\kappa}= [\kappa_{1}, \kappa_{2}, \ldots, \kappa_{Q} ]^{\top} $, $ \Omega_{M}=\{w_m \in \Omega; \ m=1,2,\ldots,M\}$ and $\mathbf{a}= [b^{+}]$\;

   Compute $\mathbf{Z} = [z_{1}, z_{2}, \ldots, z_{Q} ]^{\top}$\ according to \eqref{ZZ}\;
   Choose  sample set $\Omega_{M_j}=\{w_m \in \Omega_M; \ m=1,2,\ldots,M_j\}$ for $\kappa_j$, tolerance $\delta$ and maximum iteration $T$\;

   \For{ $j=1,2,\ldots,Q$}{
		\For{ $m=1,2,\ldots,M_j$}{
Set $a_{m,p}=0,p=2z_{j-1}+1, 2z_{j-1}+2, \ldots,2z_j$\;
		Define $\mathbf{a}_{m,z_j}= [\mathbf{a}^{\top}, a_{m,2z_{j-1}+1}, \ldots, a_{m,2z_j} ]^{\top}$\;

		Set $\mathbf{a^{(0)}}_{m,z_j}= \mathbf{a}_{m,z_j}$, $\delta_1=1$ and $t=0$\;

        \While{ $\delta_1 > \delta $ {\rm and} $t\leq T$}{
        {\bf Step 1.} \\
         Compute $A=\big[A^{(1)}; \ldots;A^{(L)}\big],
\bm{\Psi}=\big[\bm{\Psi}^{(1)}; \ldots; \bm{\Psi}^{(L)} \big],
G=\big[G^{(1)}; \ldots; G^{(L)}\big]$\;
		Solve the regularization equation $\left(A+\varepsilon I\right)\Psi=G$\;
        {\bf Step 2.} \\
        Compute  $\mathbf{J}(\mathbf{a}_{m,z_j})=[J
        ^{(1)}(\mathbf{a}_{m,z_j}),\ldots,
J^{(L)}(\mathbf{a}_{m,z_j})]^{\top}$\;
Compute the Jacobi matrix $\mathbf{D}\mathbf{J}=\left(\dfrac{\partial J
^{(l)}}{\partial a_{m,p}} \right)_{l=1,2,\ldots,L, p=0,1,\ldots,2z_j}$\;
$\mathbf{a}_{m,z_j}^{(t+1)}=\mathbf{a}_{m,z_j}^{(t)}-\eta_{\kappa_j}\mathbf{D}\mathbf{J}^{\top}(\mathbf{a}_{m,z_j}^{(t)})\mathbf{J}(\mathbf{a}_{m,z_j}^{(t)})$\;
$\delta_1=\|\mathbf{a}_{m,z_j}^{(t+1)}-\mathbf{a}_{m,z_j}^{(t)} \| $\;
$t=t+1$\;
}
$\tilde{\mathbf{a}}_{m,z_j}=\mathbf{a}^{(t)}_{m,z_j}$\;
$f_{m,z_j}(x)=\tilde{a}_{m,0}+\sum_{p=1}^{z_j}\left[\tilde{a}_{m,2p-1}\cos(px)+\tilde{a}_{m,2p}\sin(px)   \right]$\;
}
$\mathbf{a}=\frac{1}{M_j}\sum_{m=1}^{M_j}\tilde{\mathbf{a}}_{m,z_j}$\;
}
Compute mean function $\overline{f}(x)=\frac{1}{M}\sum_{m=1}^{M}f_{m,z_Q}(x)$\;
Compute covariance matrix $\bar{C}=(\bar{c}_{ij})$\;
Compute probability density function $\widehat{f}_{d}(x)=\frac{1}{Md}\sum_{i=1}^{M}K(\frac{x-x_{i}}{d})$\;
\Return $\overline{f}, \bar{C}, $ and $\widehat{f}_{d}$.
\end{algorithm}
\section{Numerical experiments}
In this section, we give several numerical experiments to verify the effectiveness of our proposed algorithm, including both Gaussian stochastic processes and non-Gaussian stochastic processes.
To test the stability of the reconstruction method, we added some random noise to the measured data, i.e.,
\begin{align*}
p^{d}(w;x,b^{+})= p^{d}(w;x,b^{+})(1+\tau \mathrm{rand}),
\end{align*}
where $\tau$ is the noise level, the near-field measurements $p^{d}(w;x,b^{+})$ are simulated by the finite element method \cite{LinLv2025AICM}, and 'rand' denotes the random numbers between $[-1,1]$ satisfying the uniform distribution.
Since each sample is independent and identically distributed, the numerical methods can be performed in parallel.

For all numerical examples below, we fix $c=5$, and choose the parameters $\lambda=1, \mu=1, \rho=1,\rho_f=1$.
Unless otherwise stated, we take the noise level
$\tau=0.5\%$, the two regularisation parameters are respectively taken as $\gamma=10^{-6}$ and $\varepsilon=0.001$.
The relaxation parameter $\eta_{\kappa}$ depending on the wavenumber is taken as $10^{-5}/ (\kappa+\kappa_{1}+\kappa_{2})^{3}$.
The truncation parameter in \eqref{3.3.2} and \eqref{3.3.3} is chosen to be $N=15$.
The number of samples $M$ required for Monte Carlo techniques should normally be at least $10^{5}$. However, it can be reduced in practice.
In numerical experiments, we always take $M=M_Q=1000$ and $M_{j}=100 \times j, j=1,2,\ldots,Q-1$ in the modified two-step Monte Carlo continuation algorithm.
\subsection{The effect of statistical parameters on random structures}
\begin{figure}[htp]
\begin{center}

  \includegraphics[width=0.42\textwidth]{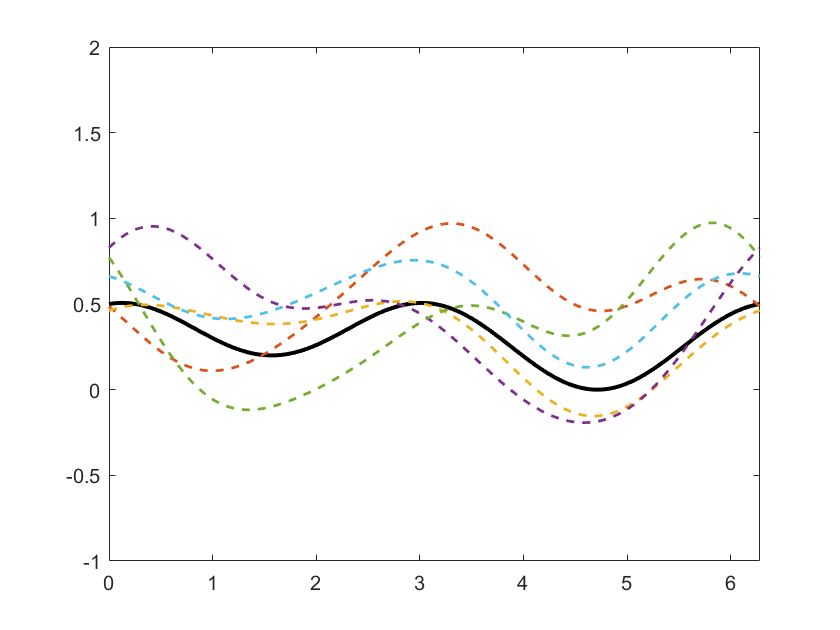}
  \includegraphics[width=0.42\textwidth]{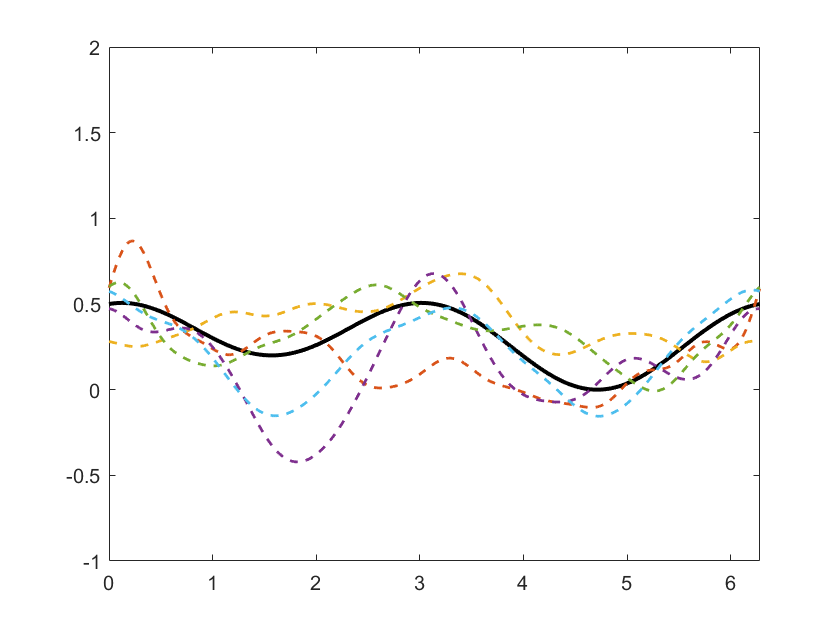}

  \includegraphics[width=0.42\textwidth]{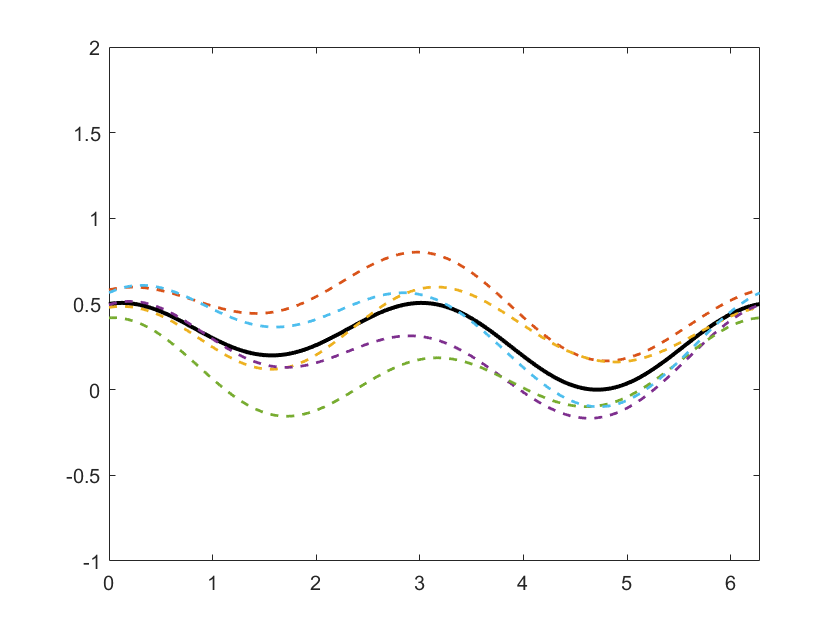}
  \includegraphics[width=0.42\textwidth]{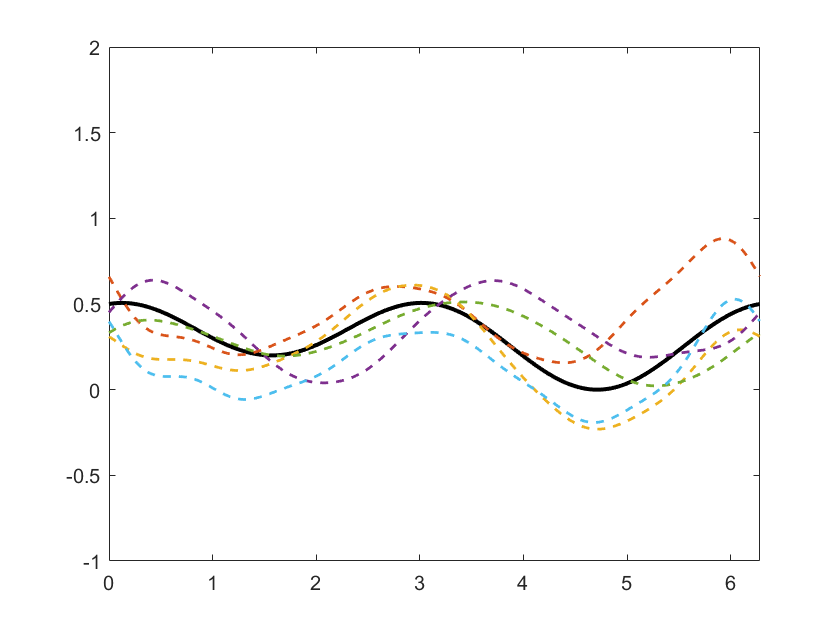}

  \includegraphics[width=0.42\textwidth]{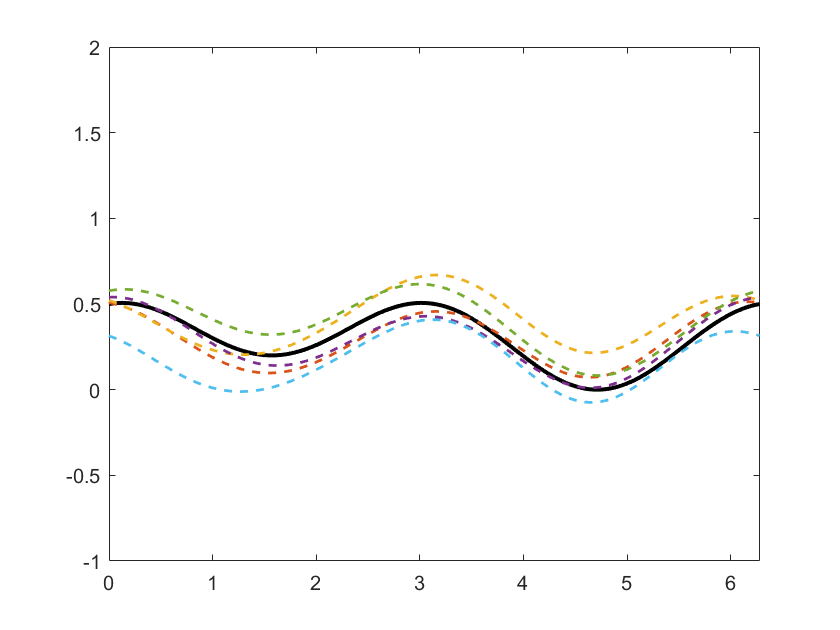}
  \includegraphics[width=0.42\textwidth]{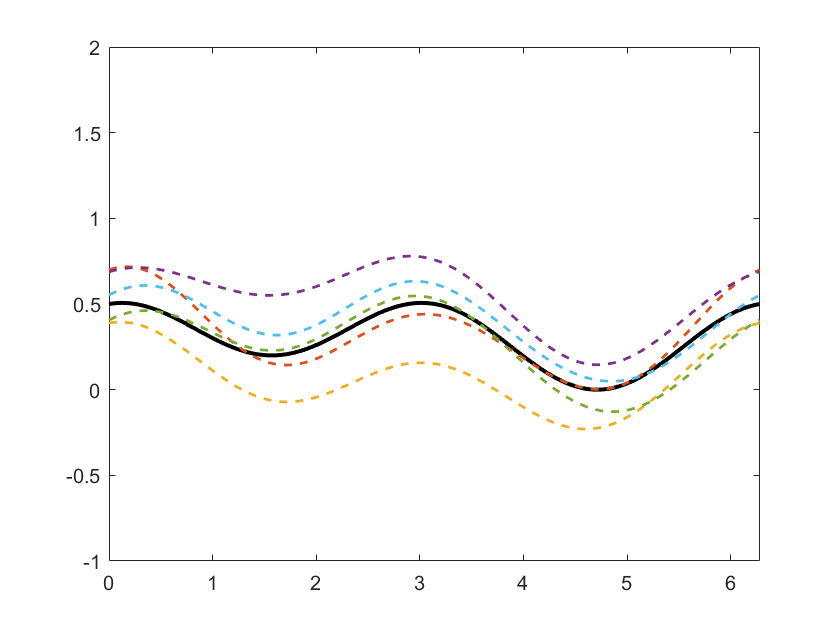}

\caption{Realisations of random surfaces in one period with different parameters.
Solid lines indicate given deterministic contour functions and dashed lines denote random samples.
Left column: $l=2$, rows 1-3 correspond to $\sigma=1/3, 1/6, 1/12$, respectively.
Right column: $\sigma=1/5$, rows 1-3 correspond to $l=0.5,1,2$, respectively.
}\label{rmsl}
\end{center}
\end{figure}

Here we discuss the effect of different statistical parameters on the random structures.
Take the covariance function
\begin{align*}
c(x-y)=\sigma^{2}e^{-\frac{|x-y|^{2}}{2 l^{2}}},
\end{align*}
where $\sigma$ is the root mean square of samples and $l$ represents the correlation length of samples.
The parameter $\sigma$ is used to measure deviation of samples from the mean, and $l$ reflects the roughness of random surfaces.
In order to show the effect of these two parameters on the random structure more clearly, we present the realizations of random surfaces in one period for different parameters with the deterministic function
\begin{align*}
\tilde{f}(x)=0.3+0.1 \sin(x)+0.2 \cos(2x),
\end{align*} as shown in Fig. \ref{rmsl}.
It can be seen that a larger $\sigma$ will increase the degree of deviation of samples from the mean, and a smaller $\l$ will produce rougher surfaces.
We define $S$ and $K$ to denote the skewness and kurtosis of the stationary non-Gaussian stochastic process, where $S$ describes the asymmetry of the distribution and $K$ denotes the sharpness of the distribution.
In this paper, we default to a standard normal distribution with a skewness of $0$ and a kurtosis of $3$.
A positive (negative) $S$ indicates that the distribution is right (left) skewed, i.e., there will be more extreme values at the right (left) end of the data. The larger the absolute value of the skewness, and the more skewed the pattern of the distribution is.
If $K$ is greater (less) than $3$, it means that the overall data distribution is steeper (flatter) than the normal distribution, with a sharp (flat) peak.

\subsection{Numerical results}
In this subsection, we introduce two key indicators for evaluating the effectiveness of the reconstruction of the stochastic process.
Define
\begin{align*}
Err_{\mathrm{mean}}=\Bigg[\sum_{i=1}^{n-1}\frac{2\pi}{n-1}(\bar{f}(x_i)-\tilde{f}(x_i))^{2}\Bigg]^{\frac{1}{2}},
\end{align*}
which can be used to characterize the error between the reconstructed mean function and the original determined function. Besides, we define
\begin{align*}
Err_{\mathrm{cov}}=\dfrac{\sum\limits_{i=1}^{n}\sum\limits_{j=1}^{n} {\big|\overline{c}(s_i,t_j)-c(s_i,t_j)\big|} }{\sum\limits_{i=1}^{n}\sum\limits_{j=1}^{n} \big|c(s_i,t_j)\big|}
\end{align*}
to characterize the relative error between the computed covariance matrix and the true covariance matrix.

\begin{Example}\label{Example1}
\textup{The first example is a stationary Gaussian stochastic process, whose original deterministic function is}
\begin{align*}
\tilde{f}(x)=0.3+0.1 \sin(x)+0.2 \cos(2x).
\end{align*}
\end{Example}

\begin{figure}[htp]
\begin{center}

  \includegraphics[width=0.32\textwidth]{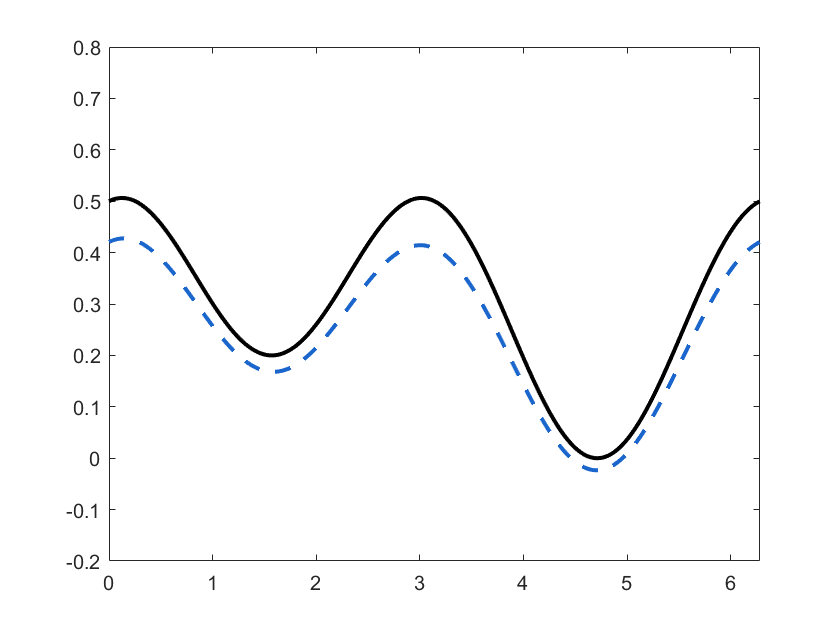}
  \includegraphics[width=0.32\textwidth]{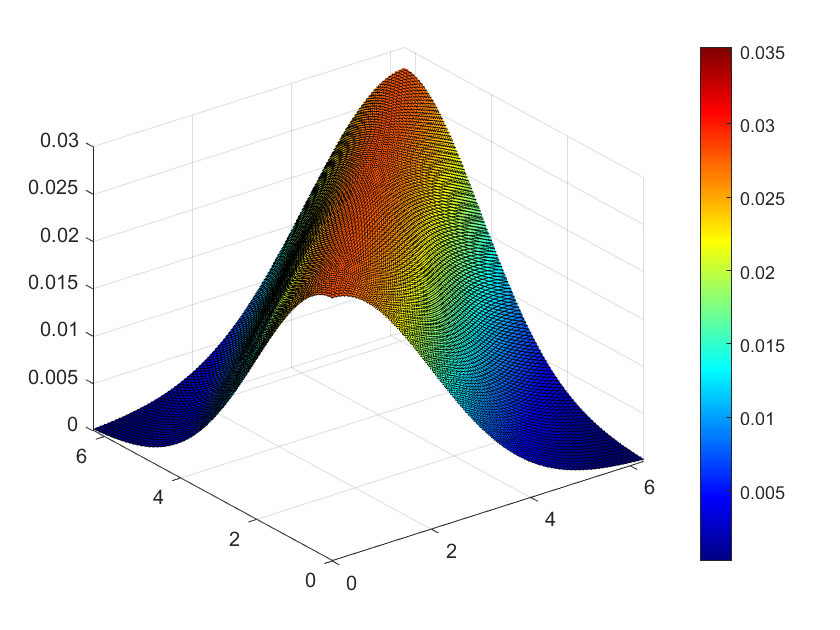}
  \includegraphics[width=0.32\textwidth]{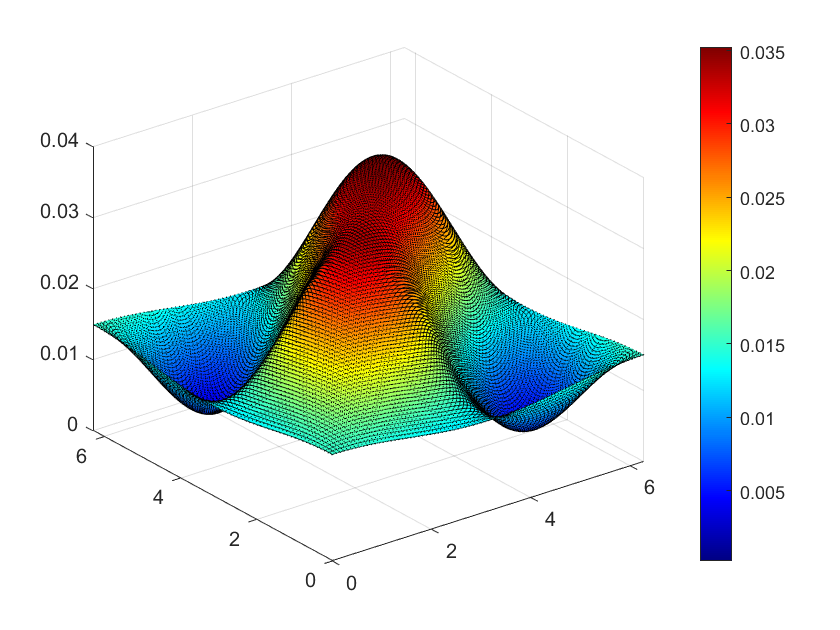}

  \includegraphics[width=0.32\textwidth]{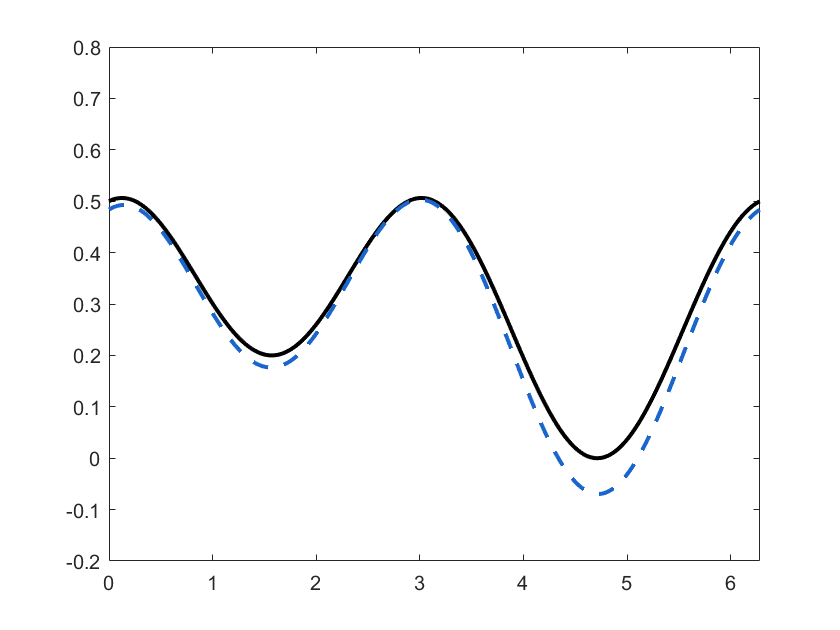}
  \includegraphics[width=0.32\textwidth]{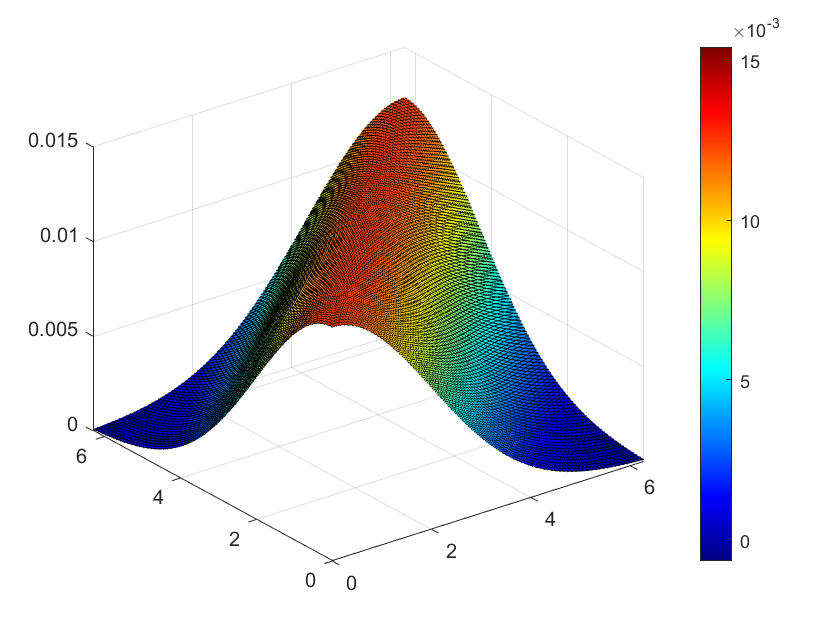}
  \includegraphics[width=0.32\textwidth]{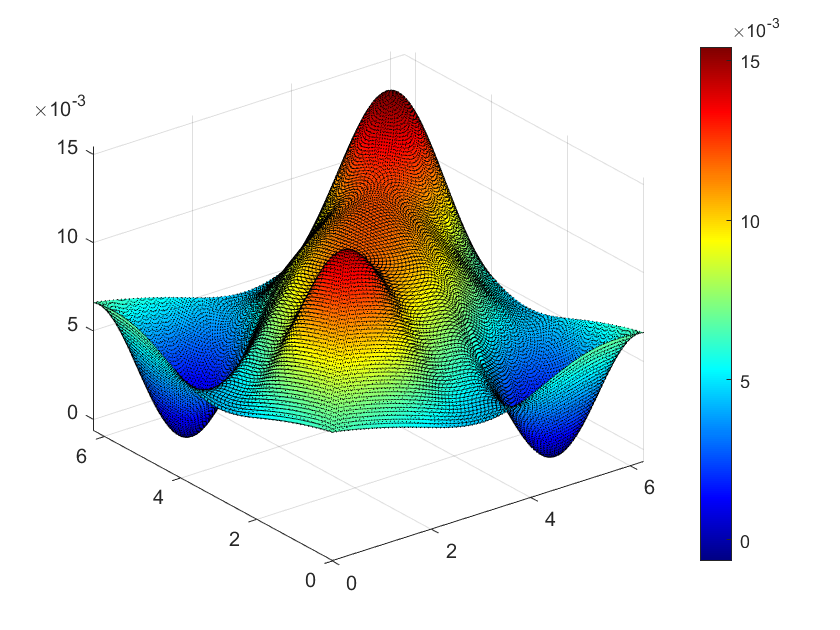}

  \includegraphics[width=0.32\textwidth]{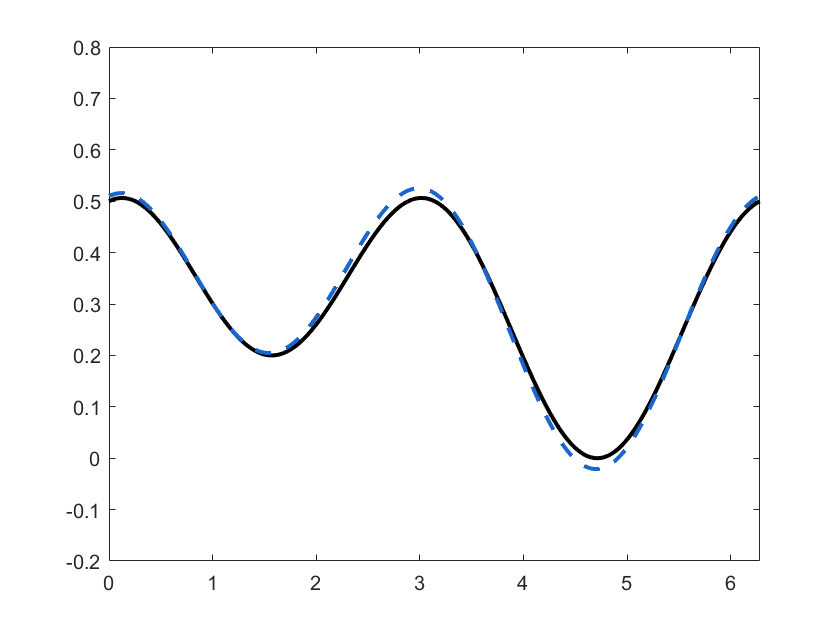}
  \includegraphics[width=0.32\textwidth]{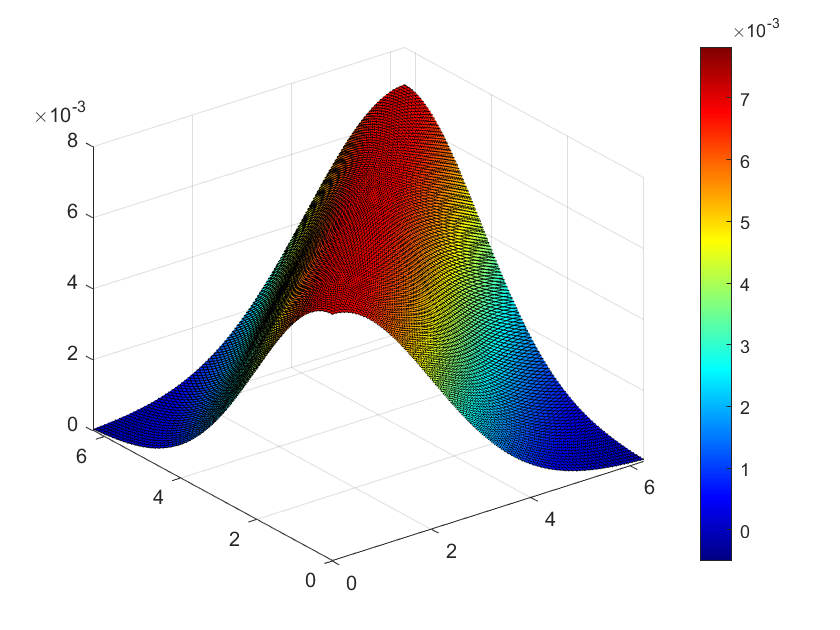}
  \includegraphics[width=0.32\textwidth]{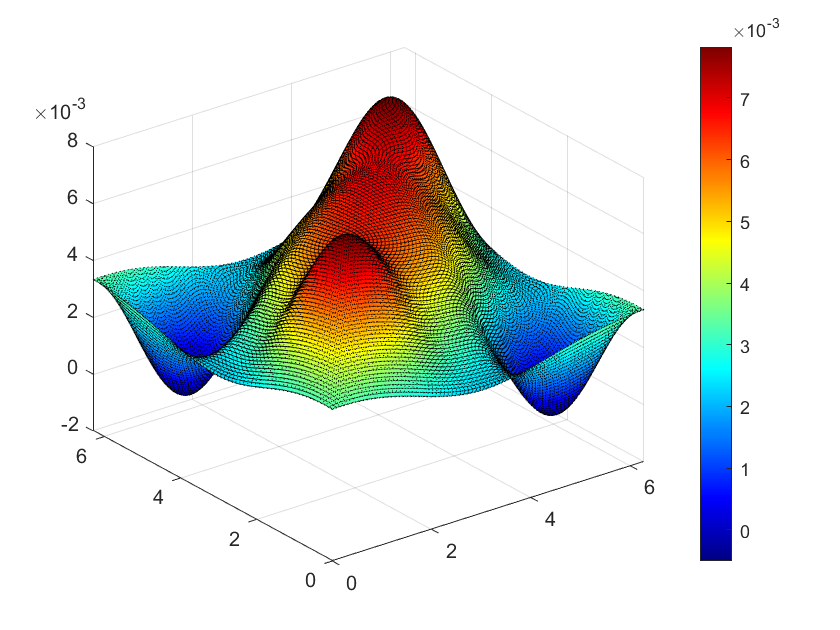}
  \caption{The reconstruction results of Example \ref{Example1}.
Left column: the solid line indicates the given deterministic profile function and the dashed line denotes the reconstructed mean profile function.
Middle column: exact covariance matrix. Right column: reconstructed covariance matrix.
Rows 1-3 correspond to $l=2$ and $\sigma=1/6, 1/9, 1/12$, respectively.
 }\label{Fig1.1}
\end{center}
\end{figure}

\begin{figure}[htp]
\begin{center}

\includegraphics[width=0.32\textwidth]{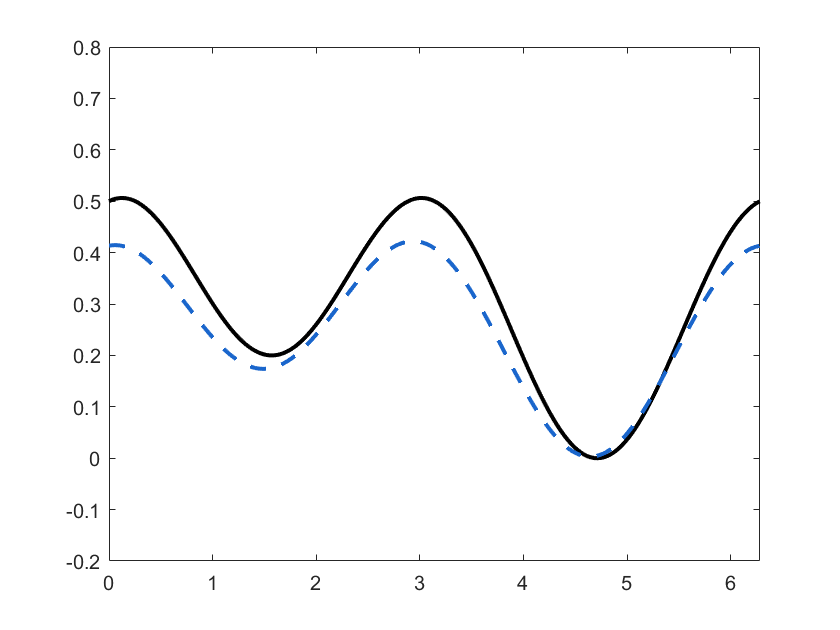}
  \includegraphics[width=0.32\textwidth]{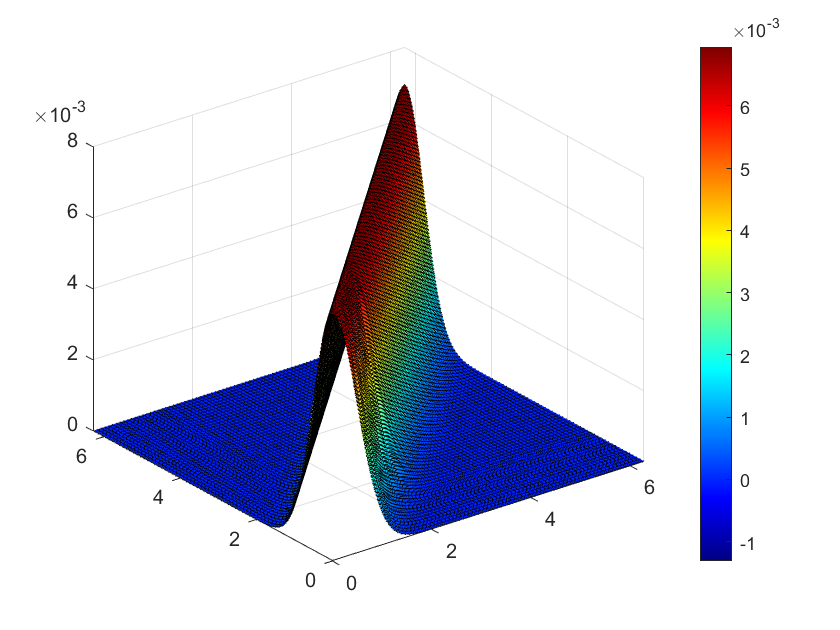}
  \includegraphics[width=0.32\textwidth]{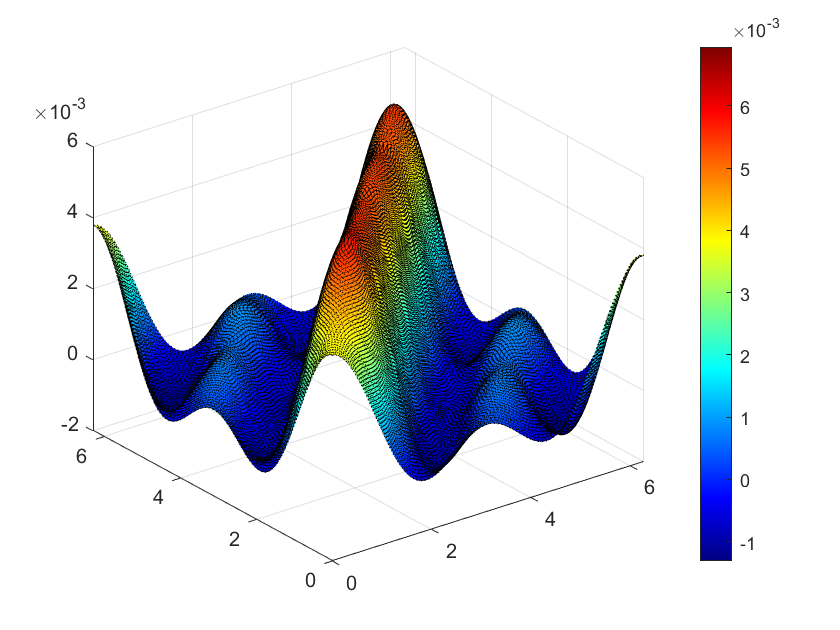}

  \includegraphics[width=0.32\textwidth]{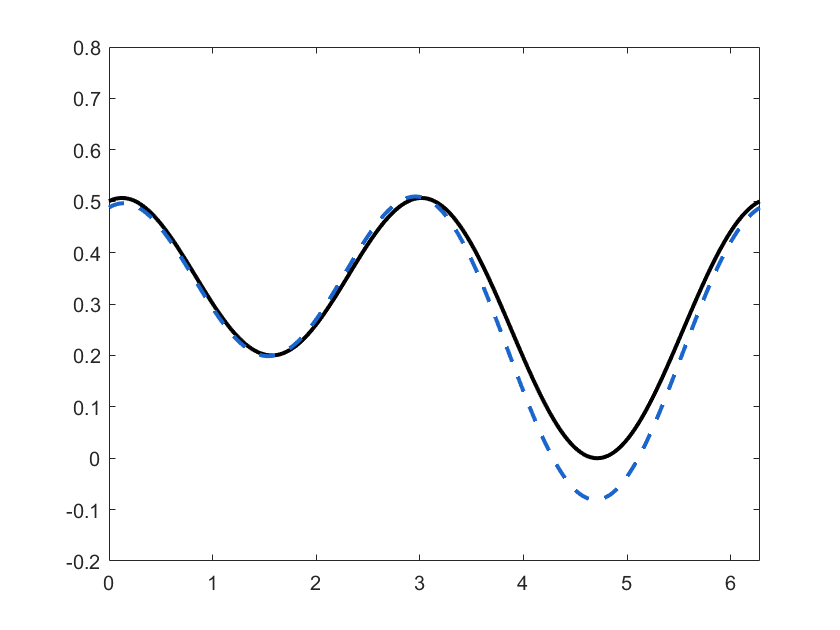}
  \includegraphics[width=0.32\textwidth]{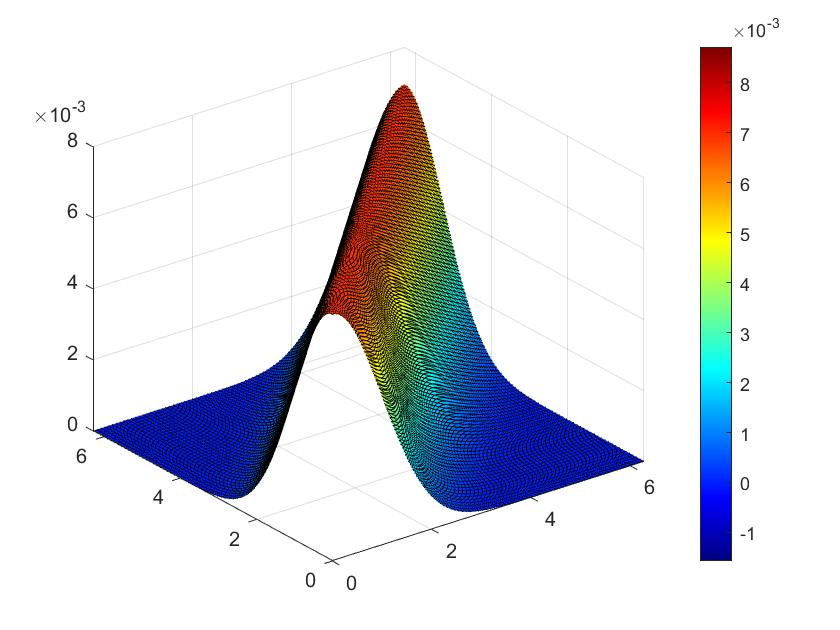}
  \includegraphics[width=0.32\textwidth]{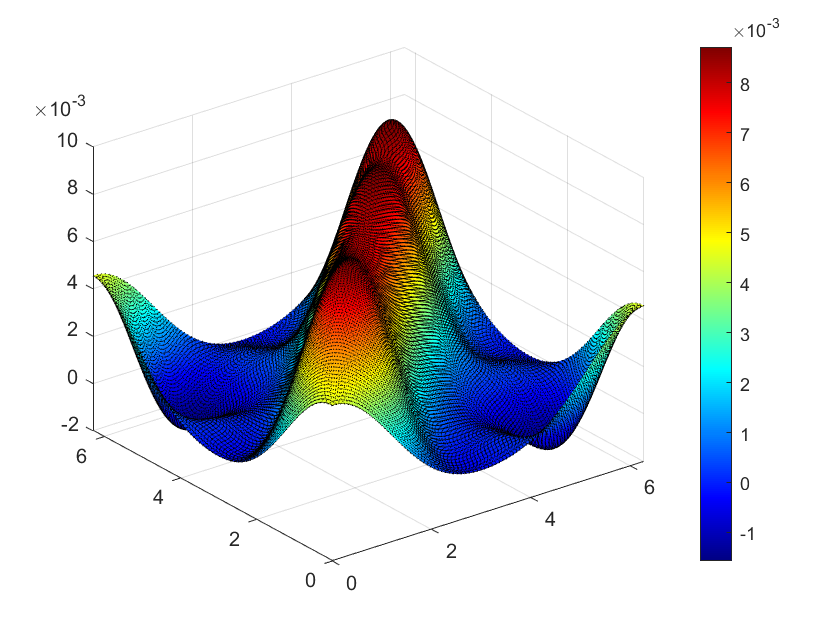}

  \includegraphics[width=0.32\textwidth]{k2_l2_s1i12_mean.png}
 \includegraphics[width=0.32\textwidth]{k2_l2_s1i12_cov1.png}
 \includegraphics[width=0.32\textwidth]{k2_l2_s1i12_cov2.png}

\caption{The reconstruction results of Example \ref{Example1}.
Left column: the solid line indicates the given deterministic profile function and the dashed line denotes the reconstructed mean profile function.
Middle column: exact covariance matrix. Right column: reconstructed covariance matrix.
Rows 1-3 correspond to  $\sigma=1/12$ and $\l=0.5,1,2$, respectively.
 }\label{Fig1.2}
\end{center}
\end{figure}

\begin{table}[htp]
\renewcommand{\arraystretch}{1.65}
\setlength{\tabcolsep}{15pt} 
\centering
\caption{The reconstruction errors of Example \ref{Example1} for $\l=2$.}\label{tab1.1}
\begin{tabular}{c|c c}
 \Xhline{1.2pt}
 $\sigma$   &$Err_{\mathrm{mean}}$ & $Err_{\mathrm{cov}}$  \\ \hline
$1/6$ & $1.52\times 10^{-1}$ & 24.24\%   \\
$1/9$ & $8.92\times 10^{-2}$ & 20.70\% \\ 
$1/12$ & $3.30\times 10^{-2}$ & 15.06\% \\  
\Xhline{1.2pt}
\end{tabular}
\end{table}

\begin{table}[htp]
\renewcommand{\arraystretch}{1.65}
\setlength{\tabcolsep}{15pt} 
\centering
\caption{The reconstruction errors of Example \ref{Example1} for $\sigma=1/12$.}\label{tab1.2}
\begin{tabular}{c|c c}
 \Xhline{1.2pt}
$\l$    &$Err_{\mathrm{mean}}$ & $Err_{\mathrm{cov}}$  \\ \hline
$0.5$ & $1.51\times 10^{-1}$ & 52.51\% \\
$1 $ & $1.00\times 10^{-1}$ & 27.87\% \\
$2$ & $3.30\times 10^{-2}$ & 15.06\% \\
\Xhline{1.2pt}
\end{tabular}
\end{table}

To investigate the effects of root mean square $\sigma$ on reconstruction results,
we take $\boldsymbol{\kappa}= [0.5, 1, 2 ]^{\top}$.
We fix $\l=2$ and take $\sigma=1/6, 1/9, 1/12$, respectively.
Fig. \ref{Fig1.1} and Table \ref{tab1.1} present the associated reconstruction results and reconstruction errors, respectively.
It can be observed that the reconstruction results perform better as $\sigma$ decreases.

Next, we explore the effects of correlation length $l$ on reconstruction results.
We set $\boldsymbol{\kappa}= [0.5, 1, 2 ]^{\top}$, fix $\sigma=1/12$ and choose $l=0.5,1,2$, respectively. At this time, the associated reconstruction results are shown in Fig. \ref{Fig1.2}.
It is apparent that a bigger $l$ will yield a better reconstruction results.
The corresponding reconstruction errors also confirm this conclusion; see Table \ref{tab1.2}.

\begin{figure}[htp]
\begin{center}

  \includegraphics[width=0.32\textwidth]{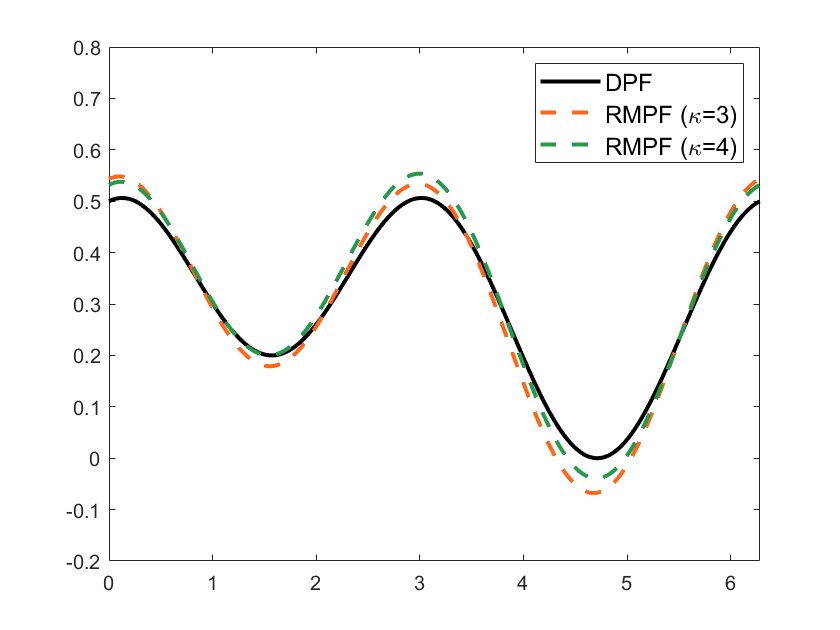}
  \includegraphics[width=0.32\textwidth]{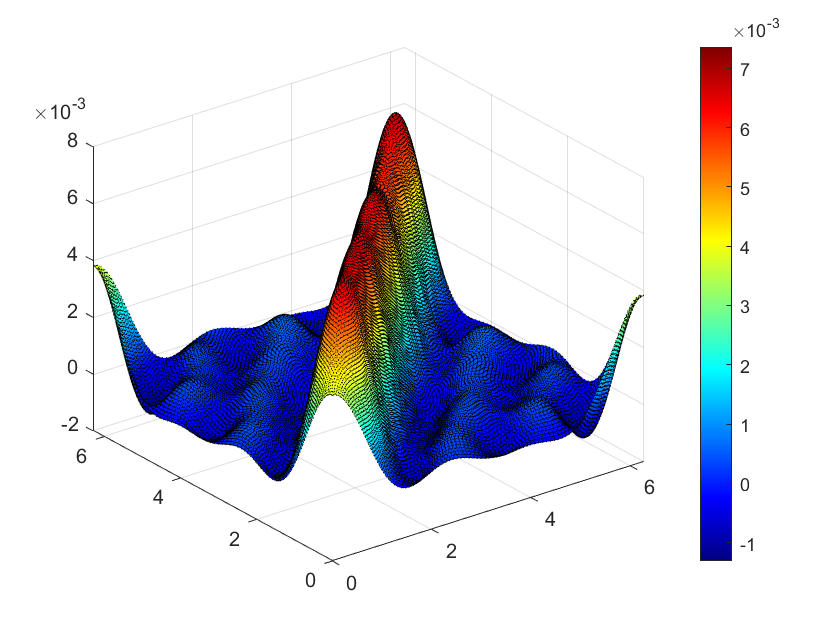}
  \includegraphics[width=0.32\textwidth]{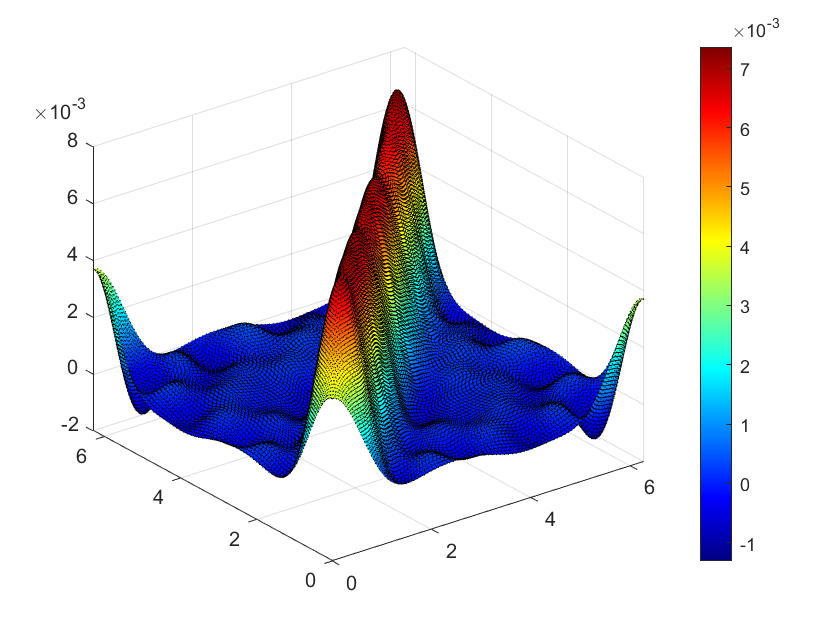}

  \caption{The reconstruction results of Example \ref{Example1}.
  Left column: the solid line indicates the given deterministic profile function (DPF) and the dashed lines denote the reconstructed mean profile functions (RMPFs).
Middle column: the reconstructed covariance matrix for $\kappa_Q=3$. Right column: the reconstructed covariance matrix for $\kappa_Q=4$.}\label{Fig1.3}
\end{center}
\end{figure}

As shown in Fig. \ref{Fig1.2} and Table \ref{tab1.2}, the reconstruction results are unsatisfactory with the parameter $l=0.5$.
To obtain a better reconstruction result, we increase the maximum of wavenumber, i.e., choose $\boldsymbol{\kappa}= [0.5, 1, 2, 3 ]^{\top}$ and $\boldsymbol{\kappa}= [0.5, 1, 2, 3, 4 ]^{\top}$.
The reconstruction results and errors are presented in Fig. \ref{Fig1.3} and Table \ref{tab1.3}, respectively.
Evidently, as the wavenumber increase, the reconstruction results are significantly improved.

\begin{table}[htp]
\renewcommand{\arraystretch}{1.65}
\setlength{\tabcolsep}{15pt} 
\centering
\caption{The reconstruction errors of Example \ref{Example1} for  $l=0.5, \sigma=1/12$.}\label{tab1.3}
\begin{tabular}{c|c c}
\Xhline{1.2pt}
$\kappa_Q$   &$Err_{\mathrm{mean}}$ & $Err_{\mathrm{cov}}$ \\ \hline
$3$ & $8.62\times 10^{-2}$ & 27.41\% \\
$4$ & $6.80\times 10^{-2}$ & 18.03\% \\
\Xhline{1.2pt}
\end{tabular}
\end{table}

\begin{Example}\label{Example2}
\textup{The second example is a static Gaussian stochastic process with an original deterministic function}
\begin{align*}
\tilde{f}(x)=0.2+0.04e^ {\cos(2x)}+0.03e^ {\cos(3x)},
\end{align*}
\textup{and take} $\sigma=1/12, \l=2$.
\textup{Obviously, the above $\tilde{f}$ contains an infinite number of Fourier modes.}
\end{Example}

\begin{figure}[htp]
\begin{center}

  \includegraphics[width=0.45\textwidth]{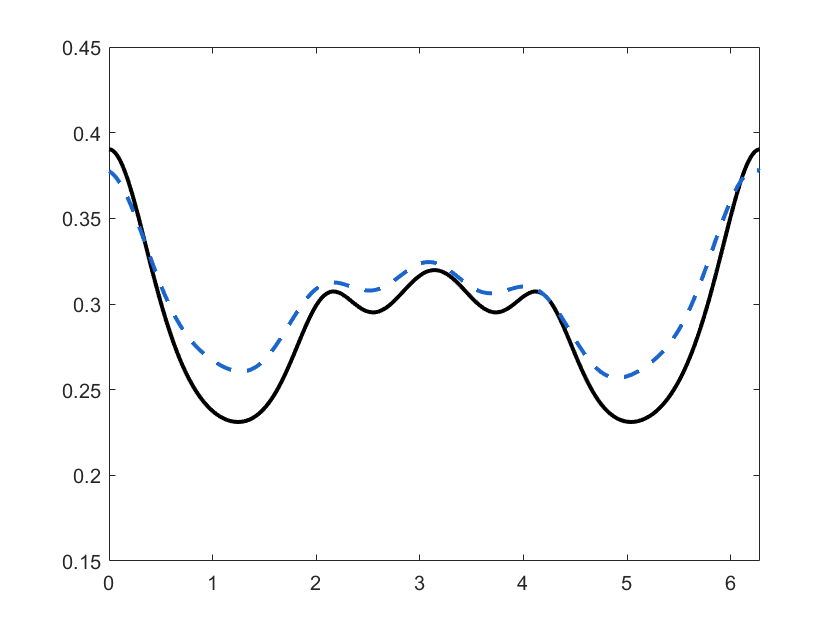}
  \includegraphics[width=0.50\textwidth]{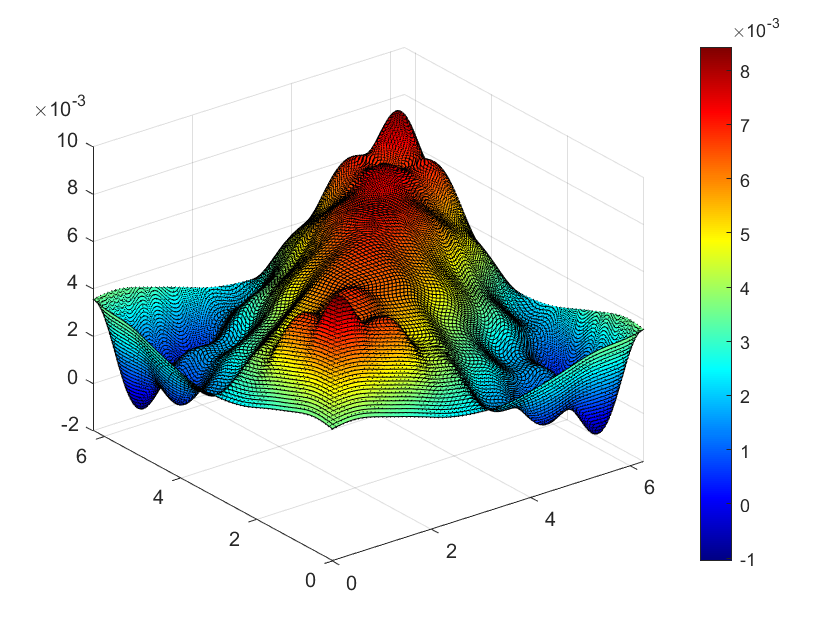}

  \includegraphics[width=0.45\textwidth]{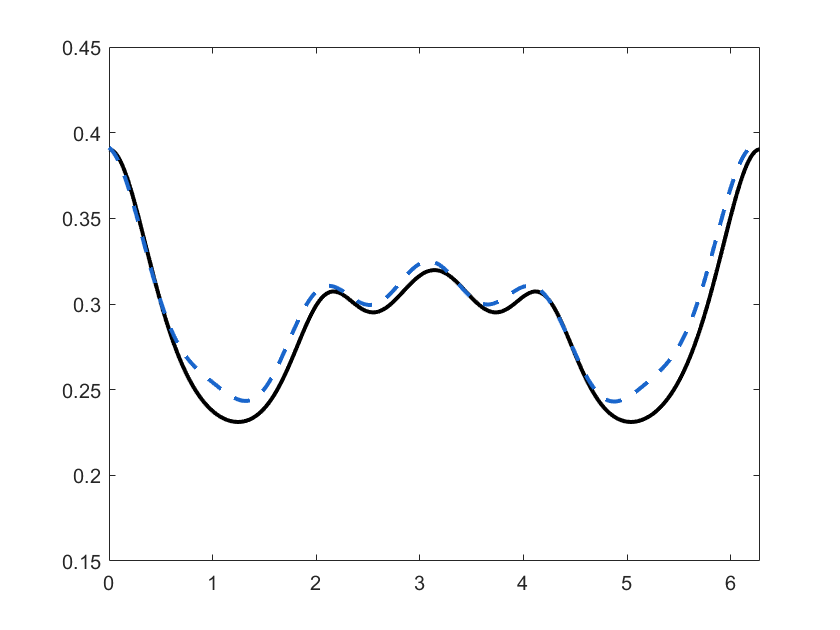}
  \includegraphics[width=0.50\textwidth]{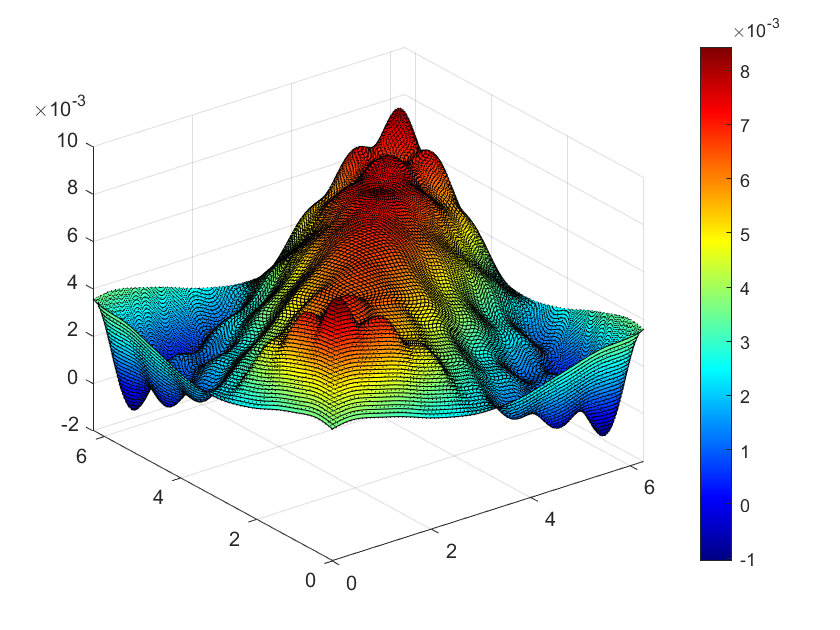}

  \includegraphics[width=0.45\textwidth]{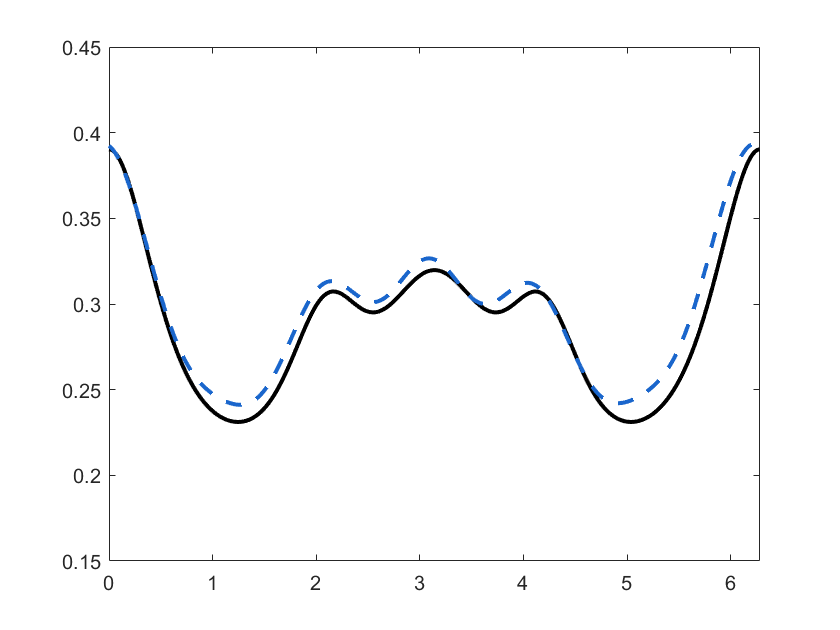}
  \includegraphics[width=0.50\textwidth]{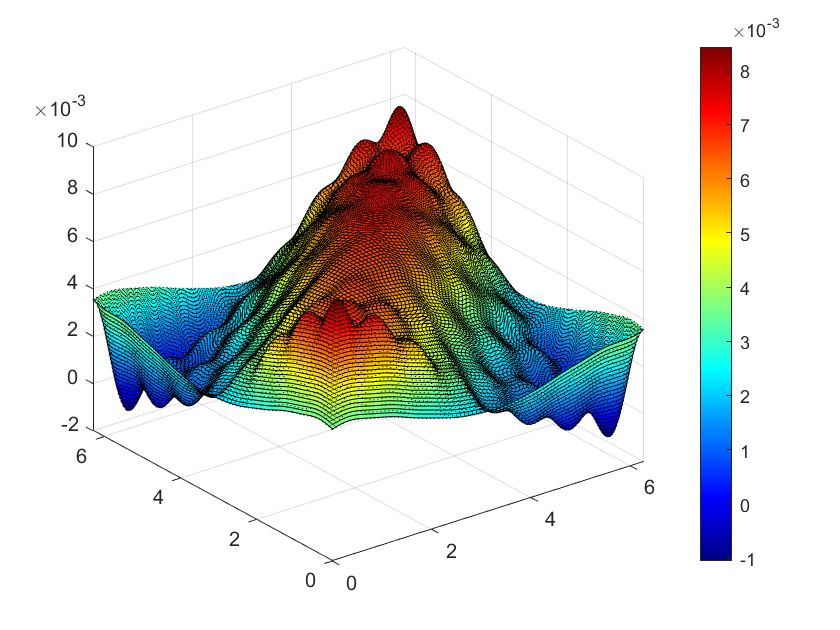}

  \caption{The reconstruction results of Example \ref{Example2}.
  Left column: the solid line indicates the given deterministic profile function and the dashed line denotes the reconstructed mean profile function.
Right column: the reconstructed covariance matrix with lines 1-3 corresponding to $\kappa_Q=6, 7, 8$, respectively.}\label{Fig2.1}
\end{center}
\end{figure}

To investigate the effects of different $\kappa_Q$ on the reconstruction, we set $\boldsymbol{\kappa}= [0. 5,1,2,3,4,5,6]^{\top}, [0. 5,1,$
$2,3,4,5,6,7]^{\top}, [0.5,1,2,$ $3,4,5,6,7,8]^{\top}$, respectively.
The associated reconstruction results  and the errors are respectively shown in Fig. \ref{Fig2.1} and Table \ref{tab2.1}.
It can be observed from Fig. \ref{Fig2.1} that all the mean profile functions and the covariance matrixes are reconstructed well.
Furthermore, a larger $\kappa_Q$ will give smaller errors of both the mean function and the covariance matrix.

\begin{table}[htp]
\renewcommand{\arraystretch}{1.65}
\setlength{\tabcolsep}{15pt} 
\centering
\caption{The reconstruction errors of Example \ref{Example2} for  $l=2, \sigma=1/12$.}\label{tab2.1}
\begin{tabular}{c|c c}
\Xhline{1.2pt}
$\kappa_Q$   &$Err_{\mathrm{mean}}$ & $Err_{\mathrm{cov}}$  \\ \hline
$6$ & $4.56\times 10^{-2}$ & 10.19\%   \\
$7$ & $2.70\times 10^{-2}$ & 9.43\% \\
$8$ & $2.62\times 10^{-2}$ & 8.82\% \\
\Xhline{1.2pt}
\end{tabular}
\end{table}


\begin{Example}\label{Example3}
\textup{
Different from the previous two cases, we consider the stationary non-Gaussian stochastic process in this example.
The original deterministic profile function is chosen as}
\begin{align*}
\tilde{f}(x)=1.2+0.1 \cos(x)+0.3\sin(2 x),
\end{align*}
\textup{and set} $S=0.9, K=5.0$.
\end{Example}

\begin{figure}[htp]
\begin{center}
  \includegraphics[width=0.50\textwidth]{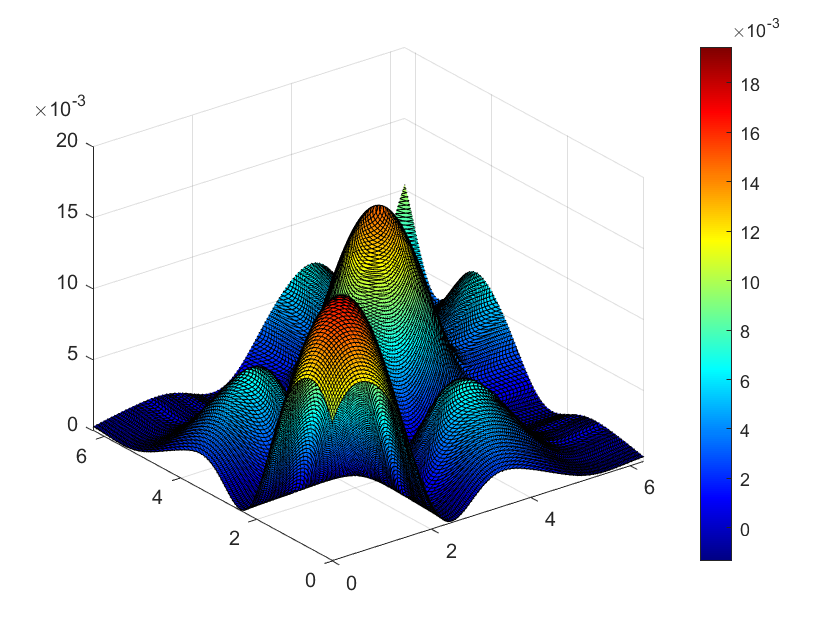}

 \caption{The exact covariance matrix of Example \ref{Example3}.}\label{Fig3.3}
\end{center}
\end{figure}

\begin{figure}[htp]
\begin{center}

  \includegraphics[width=0.45\textwidth]{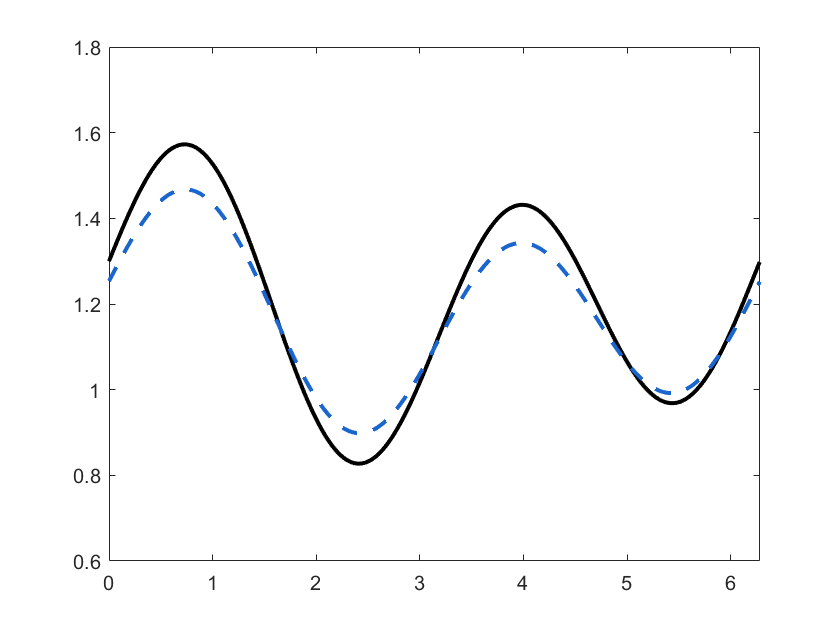}
  \includegraphics[width=0.50\textwidth]{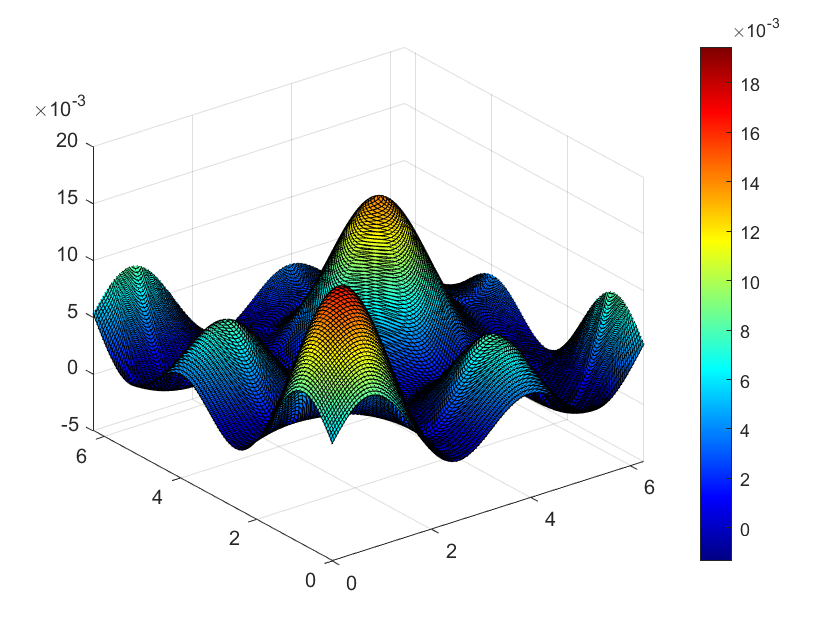}

  \includegraphics[width=0.45\textwidth]{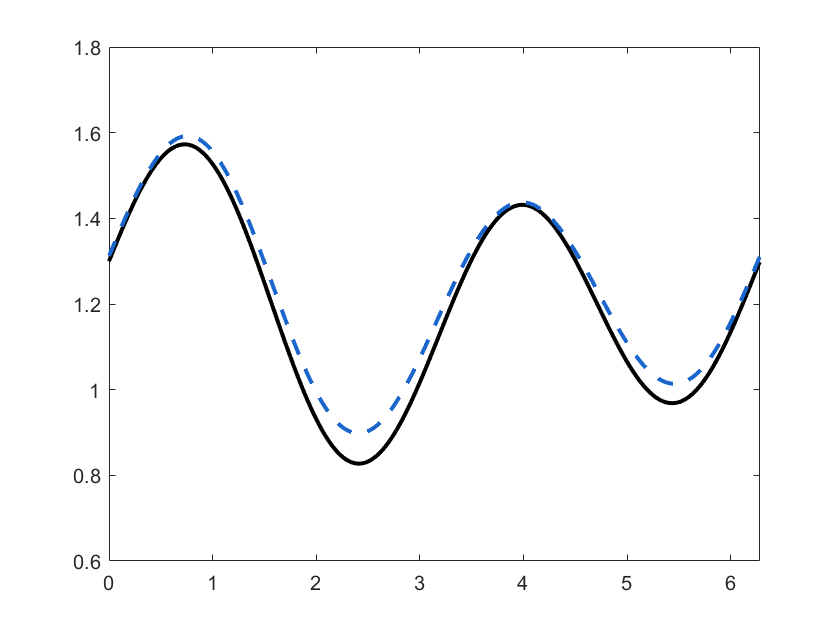}
  \includegraphics[width=0.50\textwidth]{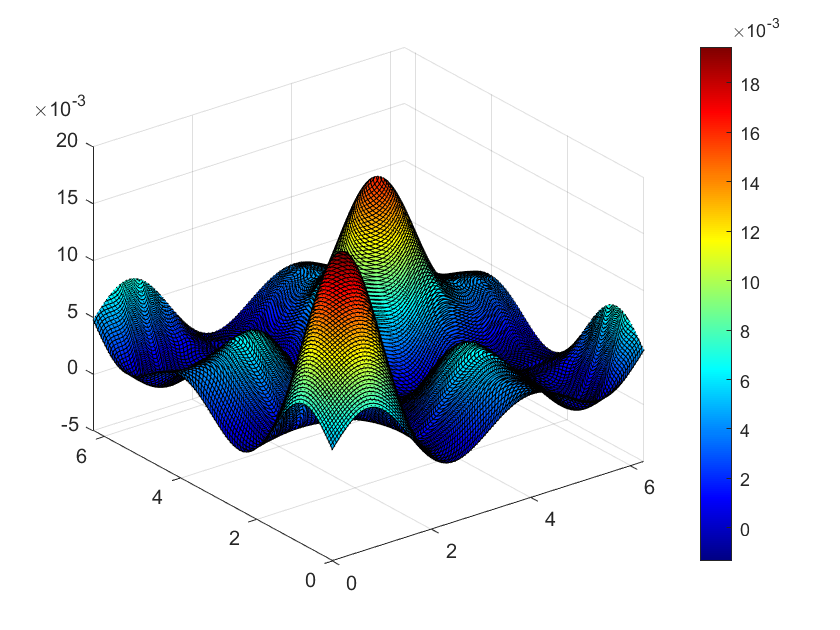}

  \includegraphics[width=0.45\textwidth]{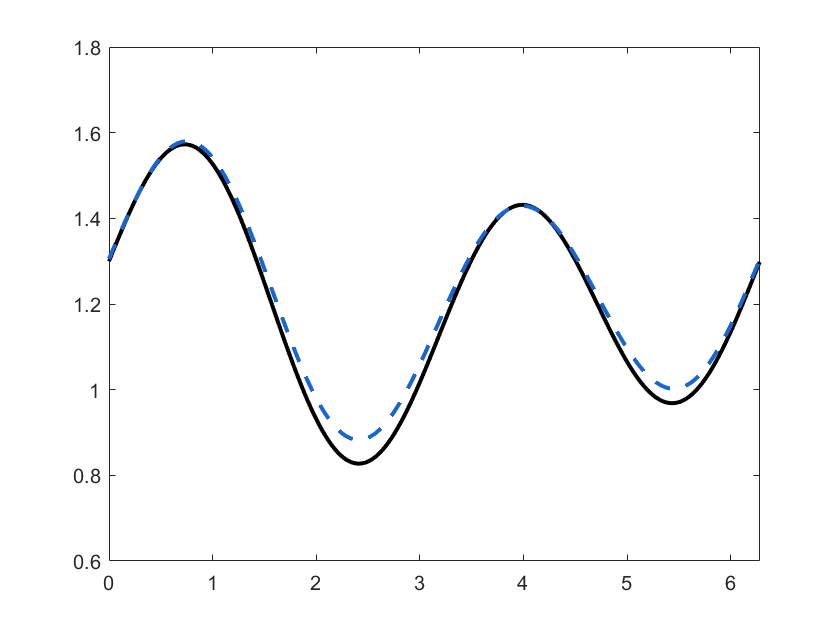}
  \includegraphics[width=0.50\textwidth]{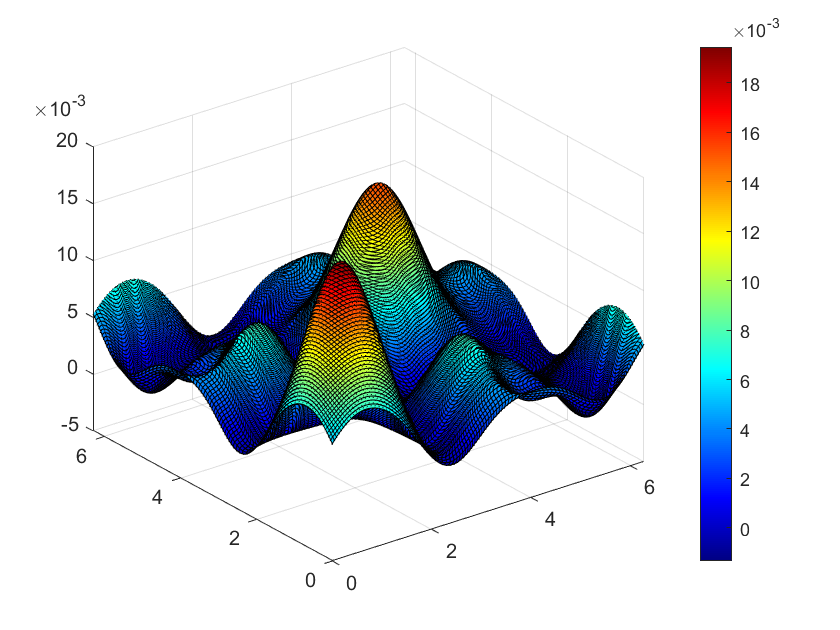}

  \caption{The reconstruction results of Example \ref{Example3}.
  Left column: the solid line indicates the given deterministic profile function and the dashed line denotes the reconstructed mean profile function.
Right column: the reconstructed covariance matrix with lines 1-3 corresponding to $\kappa_Q=2, 3, 4$, respectively.
  }\label{Fig3.1}
\end{center}
\end{figure}

\begin{figure}[htp]
\begin{center}
\subfigure[\,$x=\frac{\pi}{2}$\,]{
\includegraphics[width=0.47\textwidth]{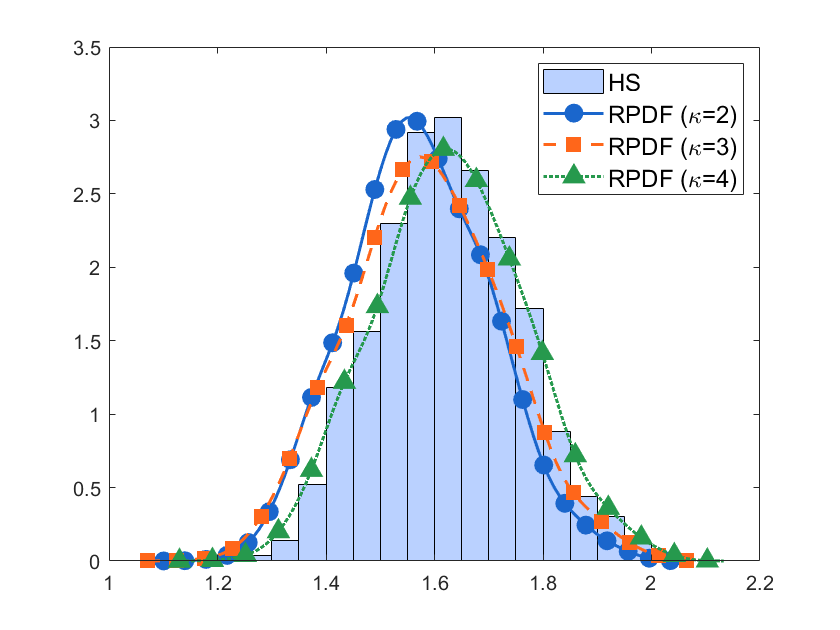}
\label{Fig3.2-1.1}
}
\subfigure[\,$x=\pi$\,]{
\includegraphics[width=0.47\textwidth]{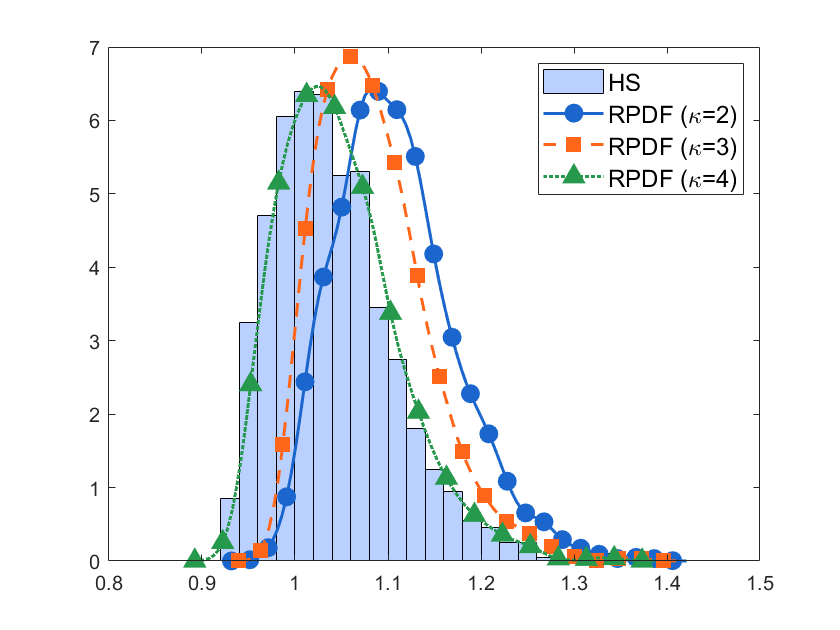}
\label{Fig3.2-1.2}
}
\subfigure[\,$x=\frac{3\pi}{2}$\,]{
\includegraphics[width=0.47\textwidth]{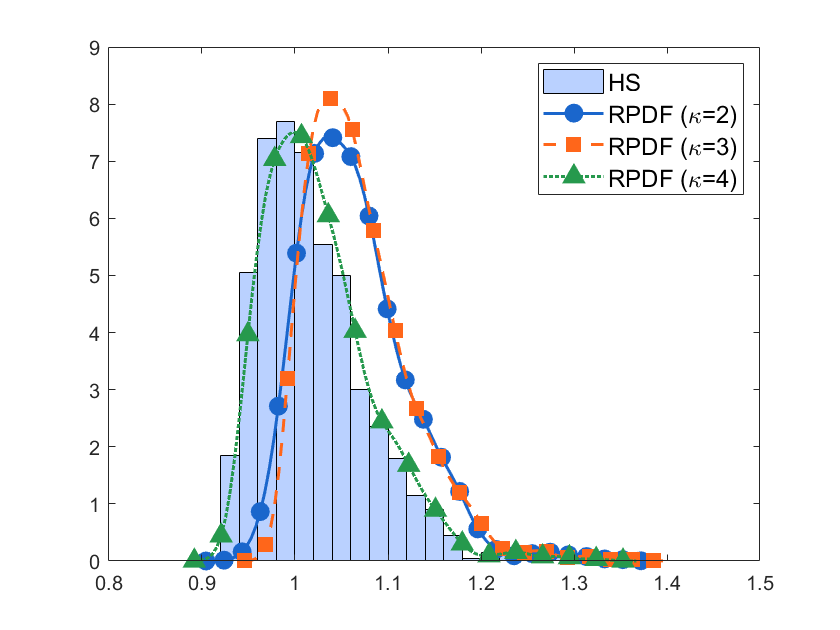}
\label{Fig3.2-1.3}
}
\subfigure[\,$x=2\pi$\,]{
\includegraphics[width=0.47\textwidth]{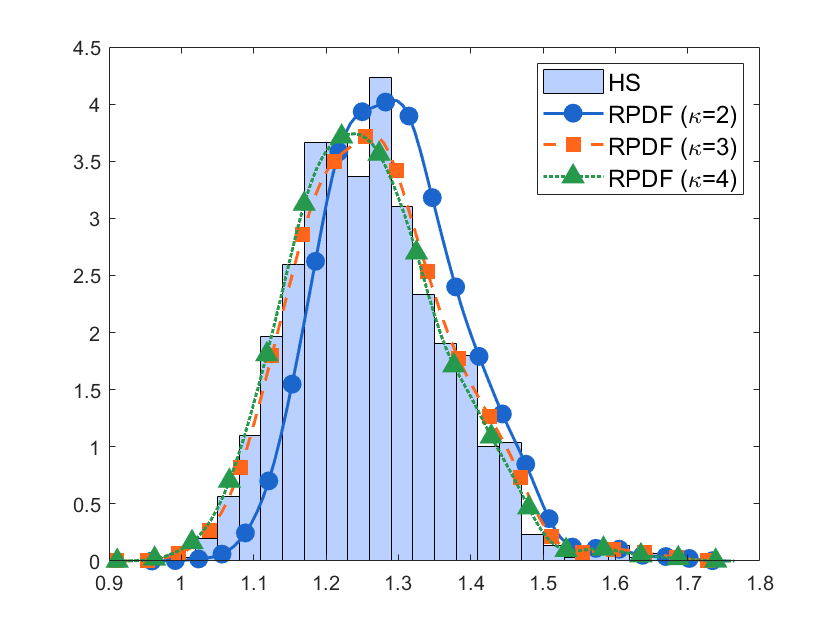}
\label{Fig3.2-1.4}
}
  \caption{Histograms of samples (HS) and results of reconstructed probability density function (RPDF) at four different positions for Example \ref{Example3}. }\label{Fig3.2}
\end{center}
\end{figure}

The exact covariance matrix is shown in Fig. \ref{Fig3.3}. Obviously,
the covariance matrix of the non-Gaussian stochastic process is more complicated than those of the previous Gaussian stochastic process.
In this case, we take $\boldsymbol{\kappa}= [0.5, 1, 2]^{\top}, [0.5, 1, 2, 3 ]^{\top}, [0.5, 1, 2, 3 ,4]^{\top}$.
The associated reconstruction results and errors are presented in Fig. \ref{Fig3.1} and Table \ref{tab3.1}, respectively.
We are ready to see that the reconstruction accuracy of the mean function and covariance matrix gets better as $\kappa_Q$ increases.
It is worth pointing out that the reconstruction results perform well enough for $\kappa_Q=4$.
In addition, we show the reconstruction of the probability density function of the stochastic process at different points in Fig. \ref{Fig3.2}. Here, we use the default kernel and optimal bandwidth in MATLAB.
For the probability density function, the approximation becomes better as $\kappa_Q$ increases.

\begin{table}[htp]
\renewcommand{\arraystretch}{1.65}
\setlength{\tabcolsep}{15pt} 
\centering
\caption{The reconstruction errors of Example \ref{Example3} for $S=0.9$ and $K=5.0$.}\label{tab3.1}
\begin{tabular}{c|c c c c}
\Xhline{1.2pt}
$\kappa_Q$   &$Err_{\mathrm{mean}}$ & $Err_{\mathrm{cov}}$ \\ \hline
$2$ & $1.48\times 10^{-1}$ & 28.86\%     \\
$3$ & $1.03\times 10^{-1}$ & 16.88\%     \\
$4$ & $7.44\times 10^{-2}$ & 13.98\%      \\
\Xhline{1.2pt}

\end{tabular}
\end{table}

\begin{Example}\label{Example4}
\textup{The fourth example is a stationary non-Gaussian stochastic process with skewness $S=1.5$ and kurtosis $K=7.0$, which are larger than those in previous example.  The original deterministic profile function is given by}
\begin{align*}
\tilde{f}(x)= 0.9-0.07e^{ \sin (x)}+0.15 e^{\cos (3x)}.
\end{align*}
\end{Example}

\begin{figure}[htp]
\begin{center}
  \includegraphics[width=0.50\textwidth]{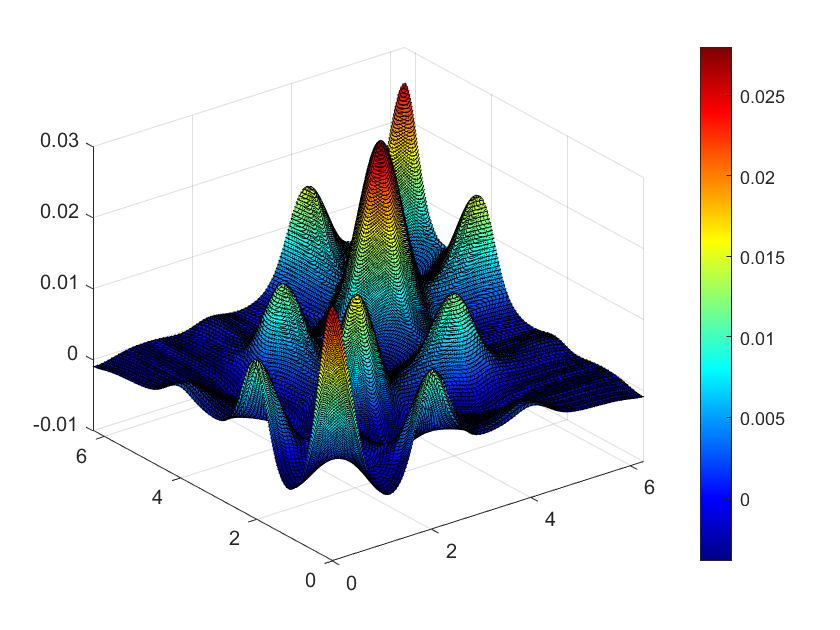}

 \caption{The exact covariance matrix of Example \ref{Example4}.}\label{Fig4.3}
\end{center}
\end{figure}

\begin{figure}
\begin{center}

  \includegraphics[width=0.45\textwidth]{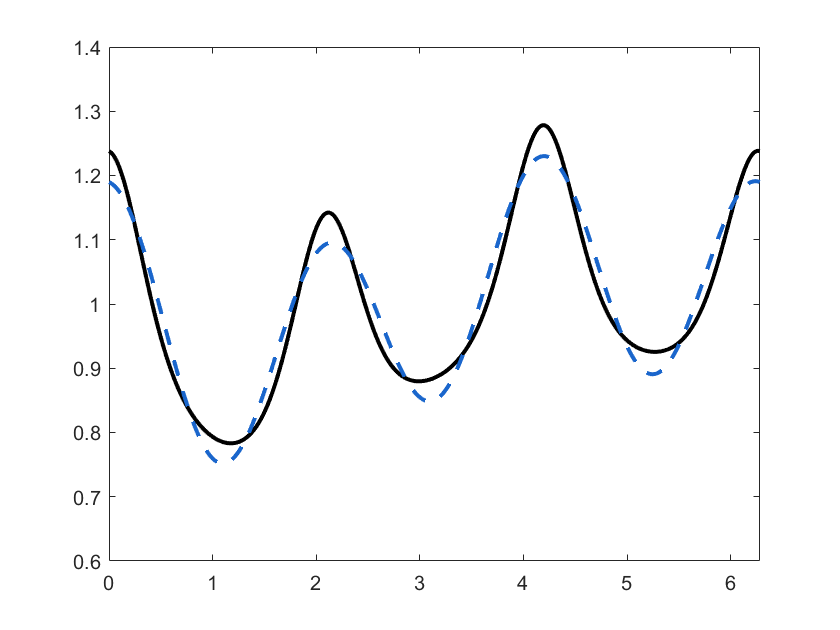}
  \includegraphics[width=0.50\textwidth]{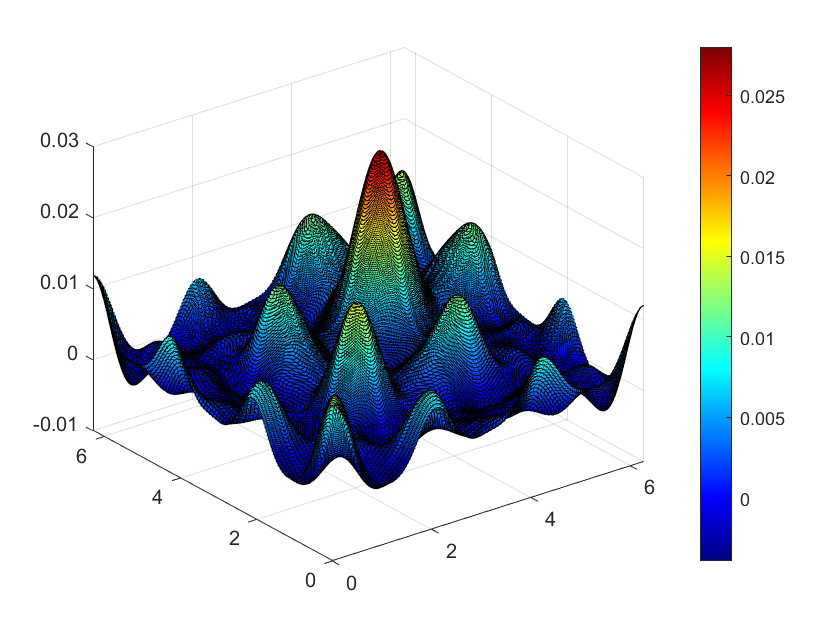}

  \includegraphics[width=0.45\textwidth]{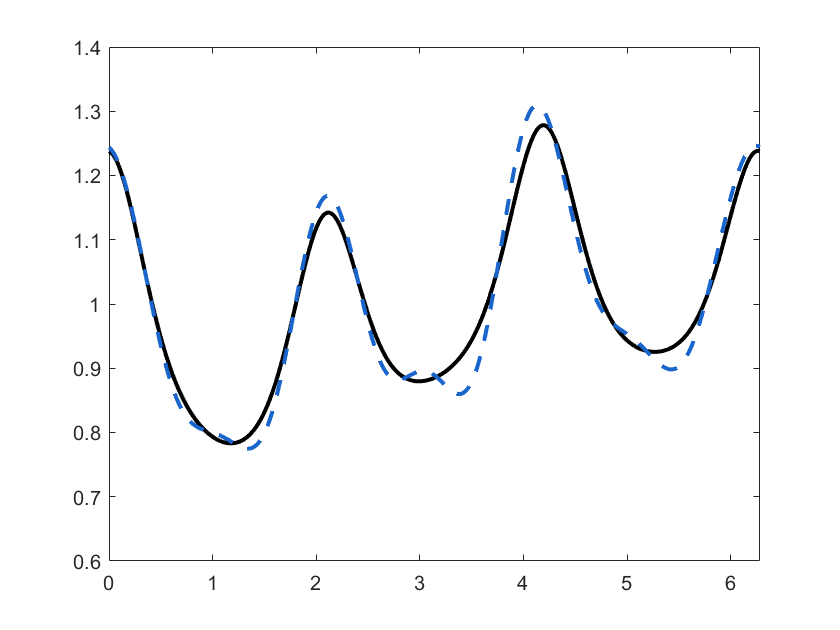}
  \includegraphics[width=0.50\textwidth]{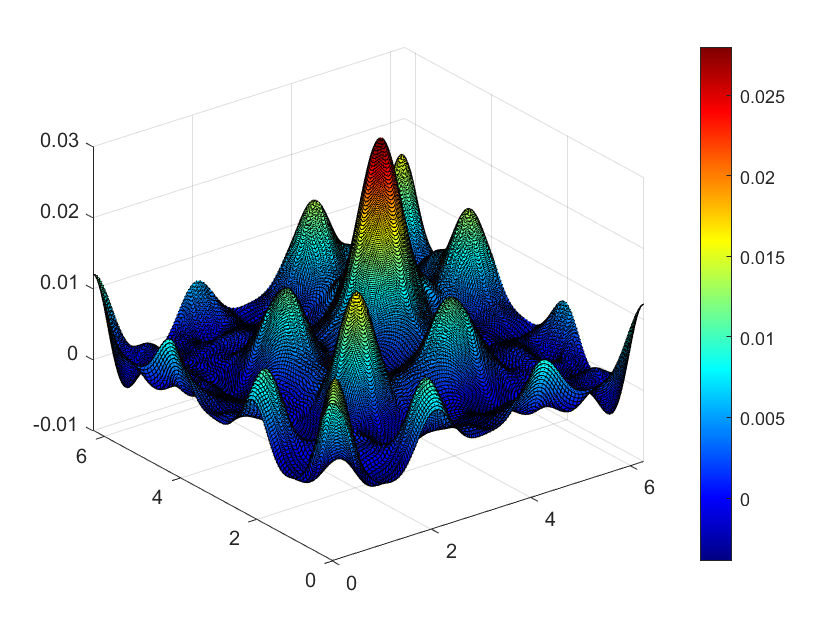}

  \includegraphics[width=0.45\textwidth]{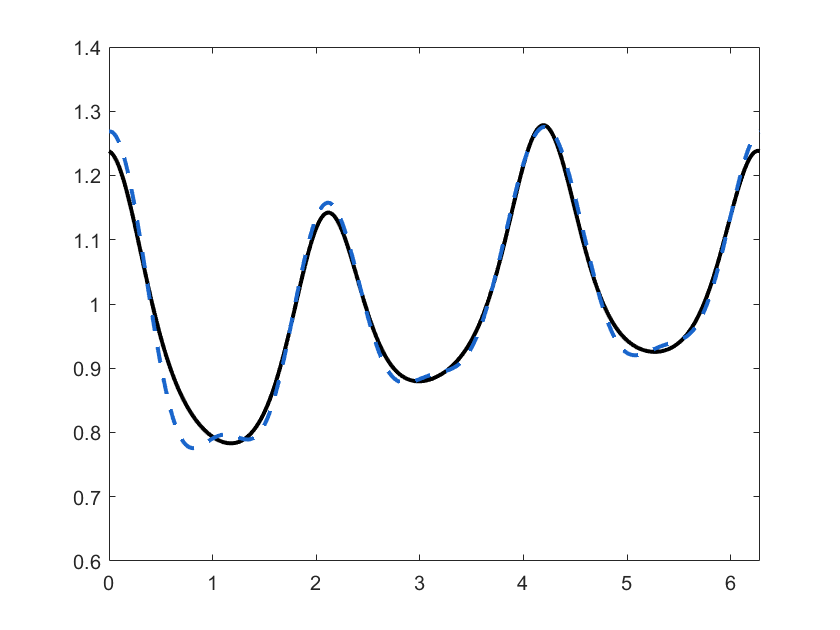}
  \includegraphics[width=0.50\textwidth]{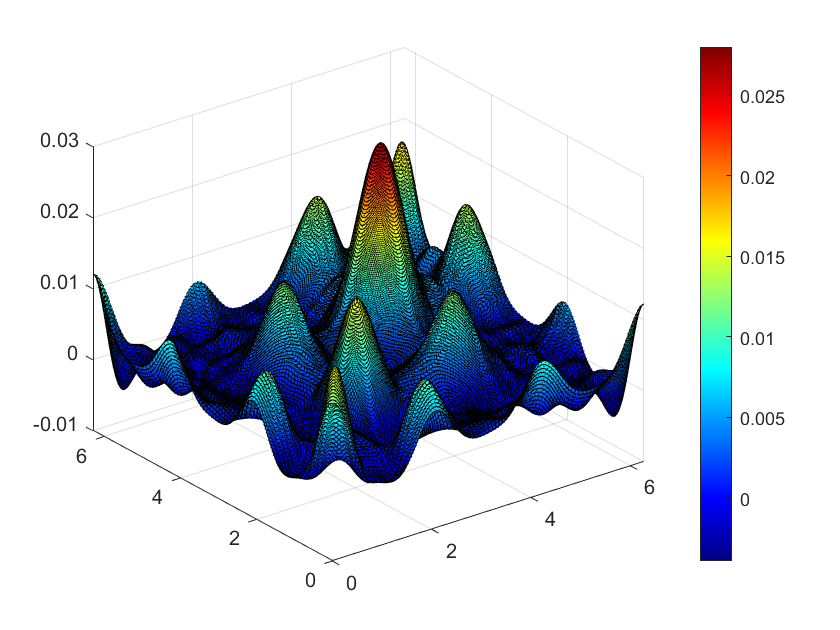}

  \caption{The reconstruction results of Example \ref{Example4}.
   Left column: the solid line indicates the given deterministic profile function and the dashed line denotes the reconstructed mean profile function.
Right column: the reconstructed covariance matrix with lines 1-3 corresponding to $\kappa_Q=6, 7, 8$, respectively. }\label{Fig4.1}
\end{center}
\end{figure}

The exact covariance matrix is plotted in Fig. \ref{Fig4.3}.
Since the deterministic function contains infinitely many Fourier modes, more iterations are required in this example.
We choose $\boldsymbol{\kappa}= [0.5, 1, 2, 3 ,4,5,6]^{\top}, [0.5,$
$1, 2, 3 ,4,5,6,7]^{\top},[0.5, 1, 2, 3 ,$
$4, 5,6,7,8]^{\top}$.
Fig. \ref{Fig4.1} and Table \ref{tab4.1} present the associated reconstruction results and errors, respectively.
It is clear to show that all three reconstruction results are satisfactory.
However, when $\kappa_Q$ is small, some deep grooves on the random surfaces cannot be reconstructed well. And the reconstructed contour function will have some burrs in the deep grooves with a too large $\kappa_Q$.
Thus, a suitable choice of $\kappa_Q$ is very important for the reconstruction.

\begin{figure}
\begin{center}
\subfigure[\,$x=\frac{\pi}{2}$\,]{
\includegraphics[width=0.47\textwidth]{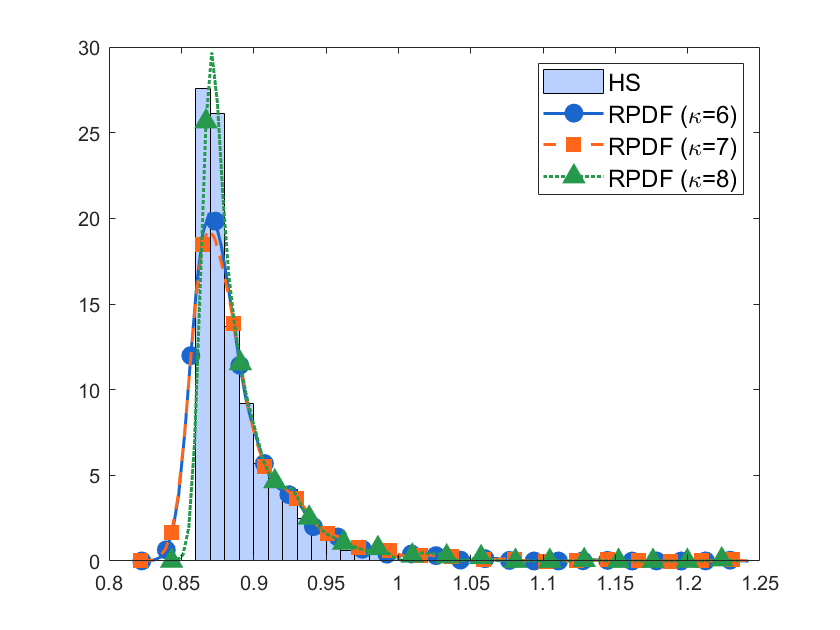}
\label{Fig4.2-1.1}
}
\subfigure[\,$x=\pi$\,]{
\includegraphics[width=0.47\textwidth]{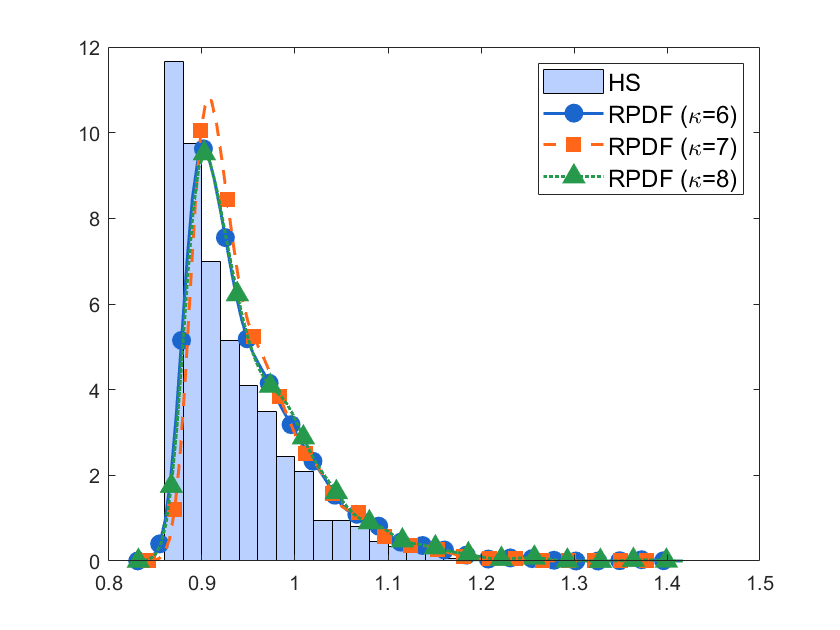}
\label{Fig4.2-1.2}
}
\subfigure[\,$x=\frac{3\pi}{2}$\,]{
\includegraphics[width=0.47\textwidth]{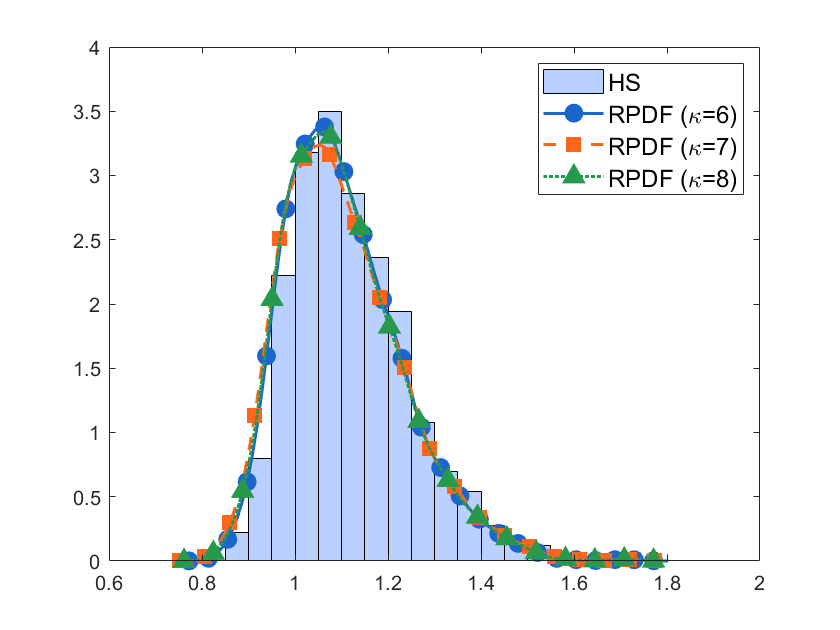}
\label{Fig4.2-1.3}
}
\subfigure[\,$x=2\pi$\,]{
\includegraphics[width=0.47\textwidth]{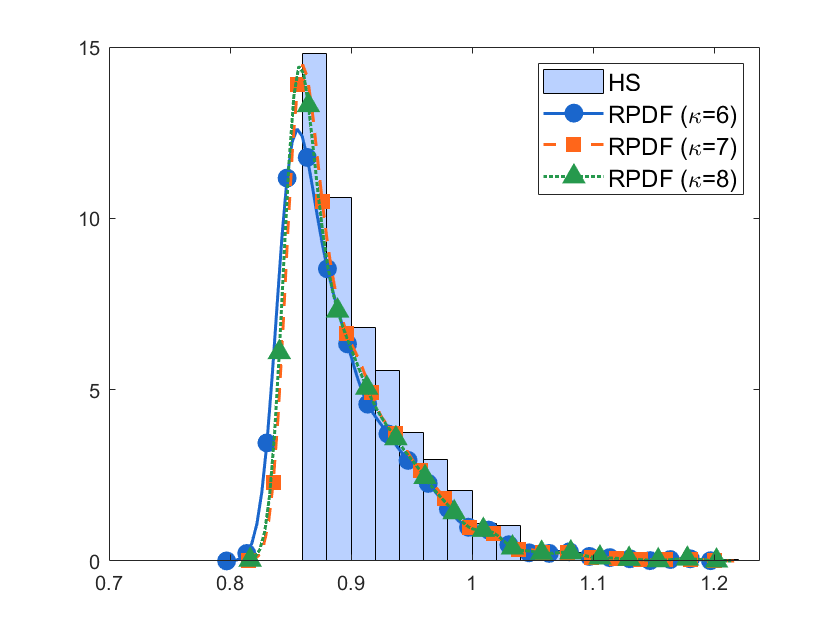}
\label{Fig4.2-1.4}
}
  \caption{Histograms of samples and results of reconstructed probability density function at four different positions for Example \ref{Example4}. }\label{Fig4.2}
\end{center}
\end{figure}

\begin{table}
\renewcommand{\arraystretch}{1.65}
\setlength{\tabcolsep}{15pt} 
\centering
\caption{The reconstruction errors of Example \ref{Example4} for $S=1.5$ and $K=7.0$.}\label{tab4.1}
\begin{tabular}{c|c c c c}
\Xhline{1.2pt}
$\kappa_Q$   &$Err_{\mathrm{mean}}$ & $Err_{\mathrm{cov}}$ \\ \hline
$6$ & $7.32\times 10^{-2}$ & 22.48\%      \\
$7$ & $6.10\times 10^{-2}$ & 17.24\%     \\
$8$ & $5.07\times 10^{-2}$ & 15.47\%     \\
\Xhline{1.2pt}

\end{tabular}
\end{table}

\begin{Example}\label{Example5}
\textup{
The last example is concerned with the binary grating which is a nonsmooth case.
We still consider a stationary non-Gaussian stochastic process with $S=0.3$ and $K=4.0$. The original deterministic function is}
\begin{align*}
\tilde{f}(x) =
\begin{cases}
1.6, &\  \mathrm{for} \  x\in [0,2]\cup [5,2 \pi] ,\\
\\
1.2, &\  \mathrm{for} \ x \in (2, 5 ).\\
\end{cases}
\end{align*}
\end{Example}

\begin{figure}
\begin{center}
  \includegraphics[width=0.50\textwidth]{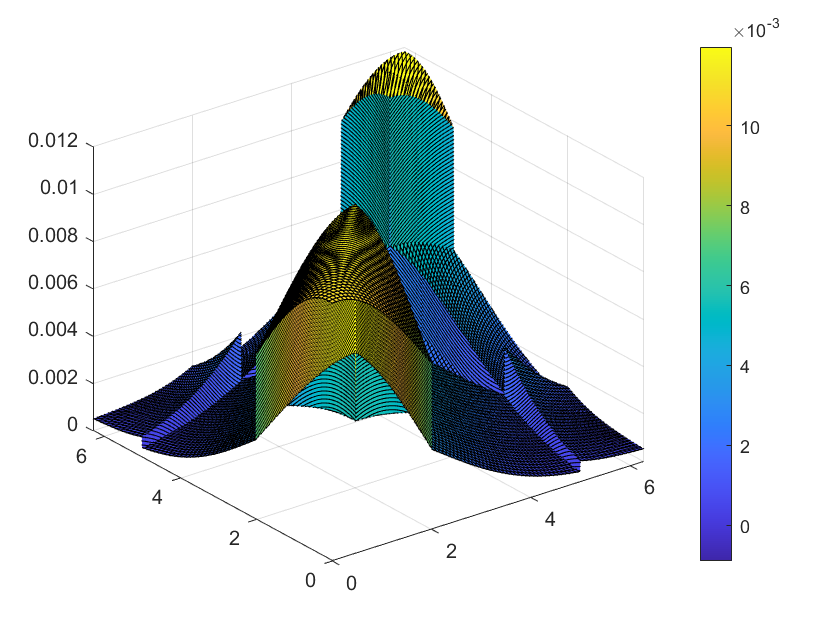}

 \caption{The exact covariance matrix of Example \ref{Example5}.}\label{Fig5.3}
\end{center}
\end{figure}

\begin{figure}
\begin{center}

  \includegraphics[width=0.45\textwidth]{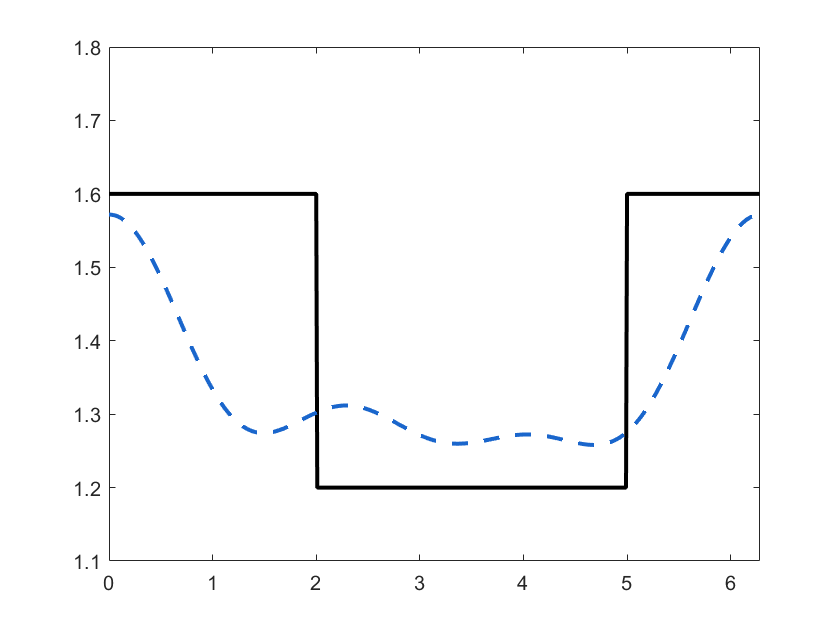}
  \includegraphics[width=0.50\textwidth]{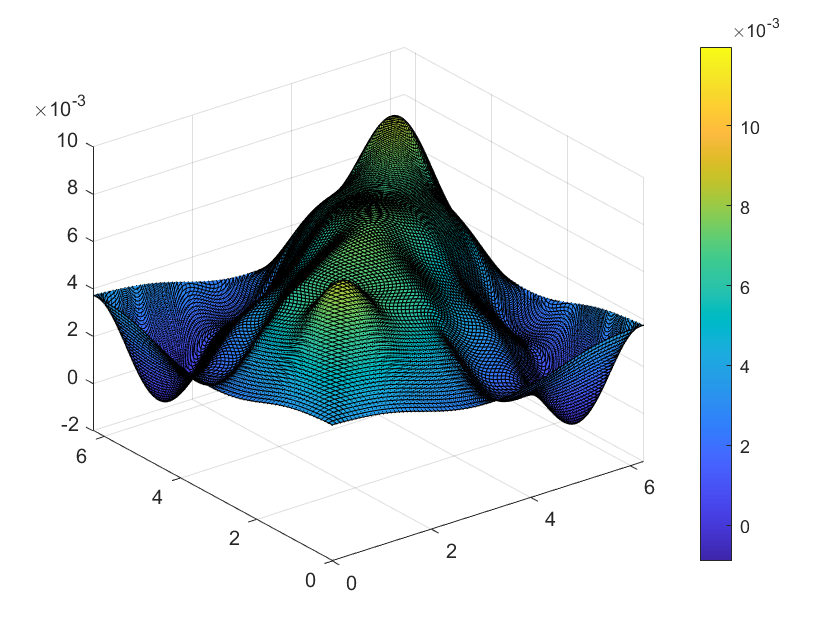}

  \includegraphics[width=0.45\textwidth]{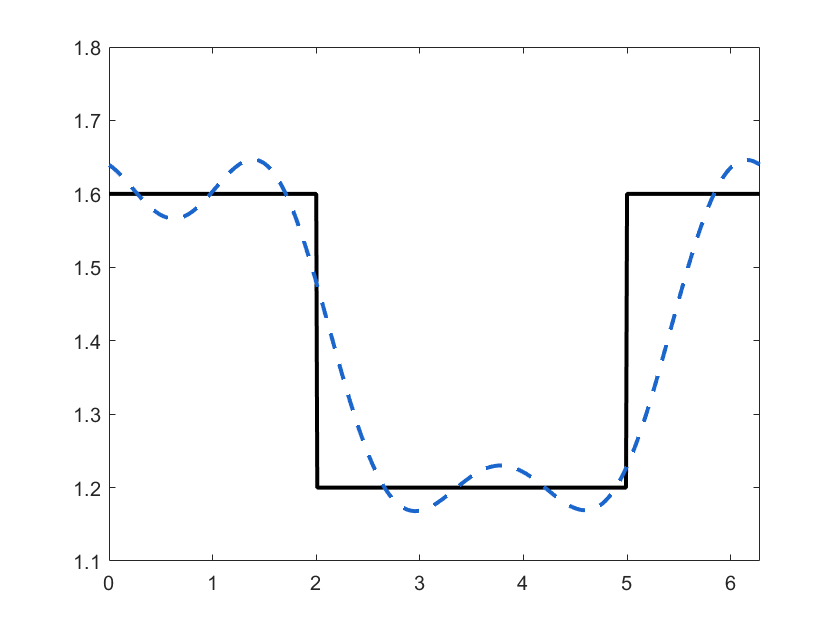}
  \includegraphics[width=0.50\textwidth]{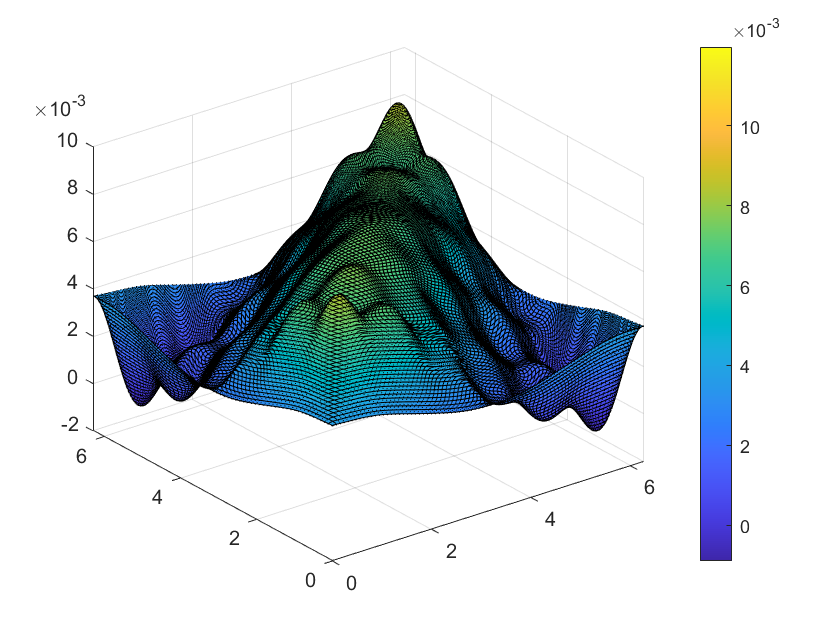}

  \includegraphics[width=0.45\textwidth]{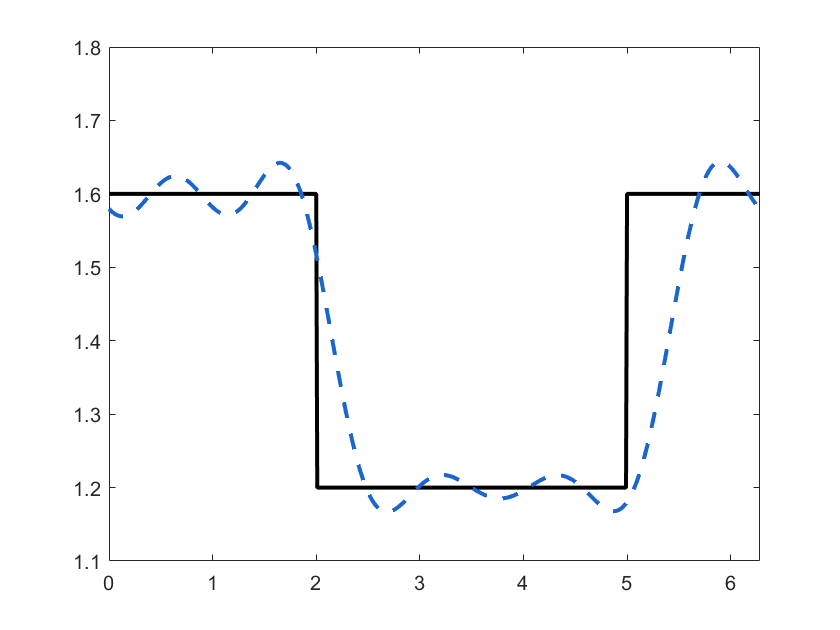}
  \includegraphics[width=0.50\textwidth]{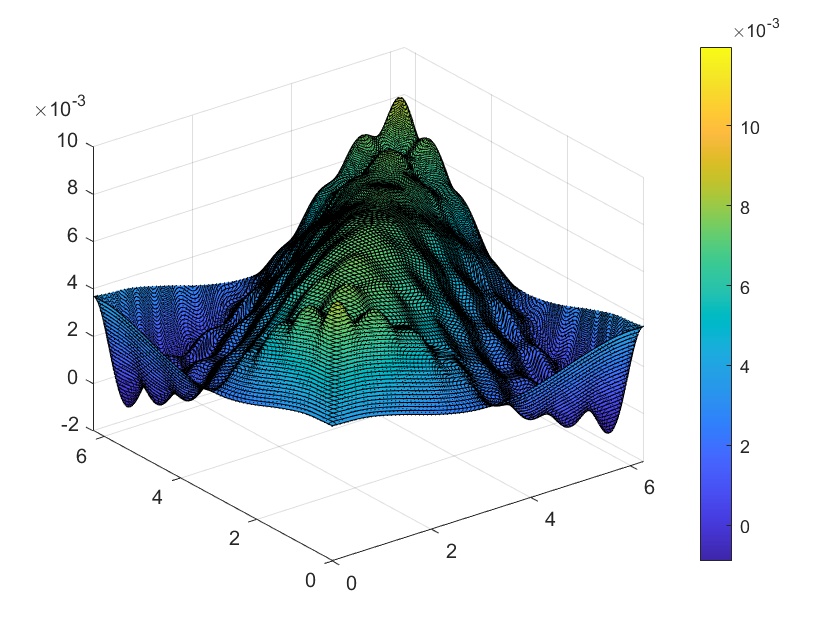}

  \caption{The reconstruction results of Example \ref{Example5}.
  Left column: The solid red line and the dashed blue line indicate the given deterministic profile function and the reconstructed mean profile function, respectively.
Middle column: exact covariance matrix. Right column: reconstructed covariance matrix. Rows 1-3 correspond to $\kappa^{+}_Q=3,5,7$, respectively. }\label{Fig5.1}
\end{center}
\end{figure}

\begin{figure}
\begin{center}
\subfigure[\,$x=\frac{\pi}{2}$\,]{
\includegraphics[width=0.47\textwidth]{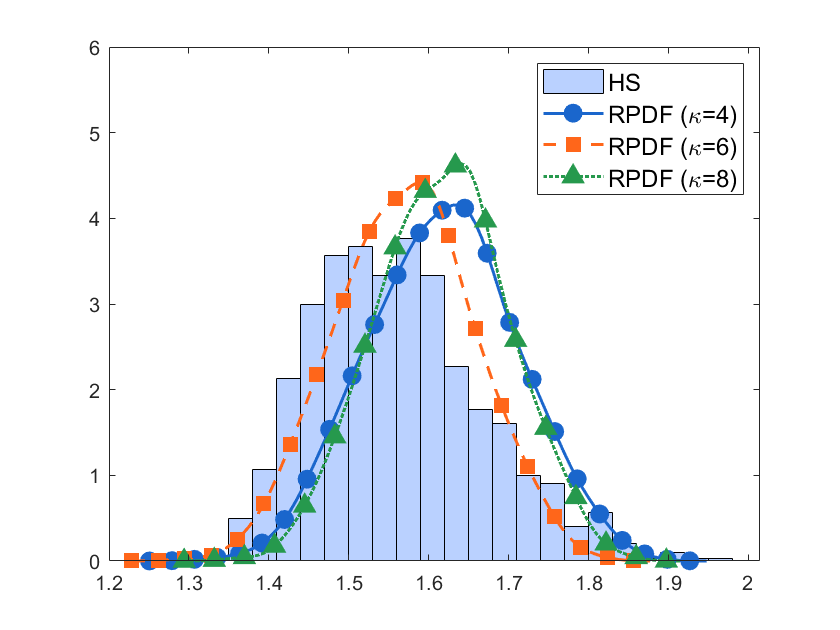}
\label{Fig4.2-1.1}
}
\subfigure[\,$x=\pi$\,]{
\includegraphics[width=0.47\textwidth]{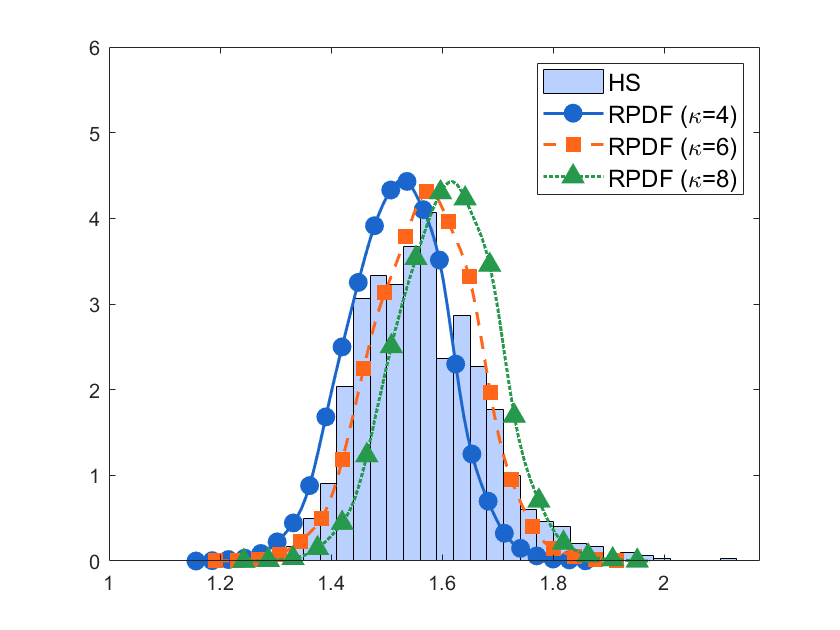}
\label{Fig4.2-1.2}
}
\subfigure[\,$x=\frac{3\pi}{2}$\,]{
\includegraphics[width=0.47\textwidth]{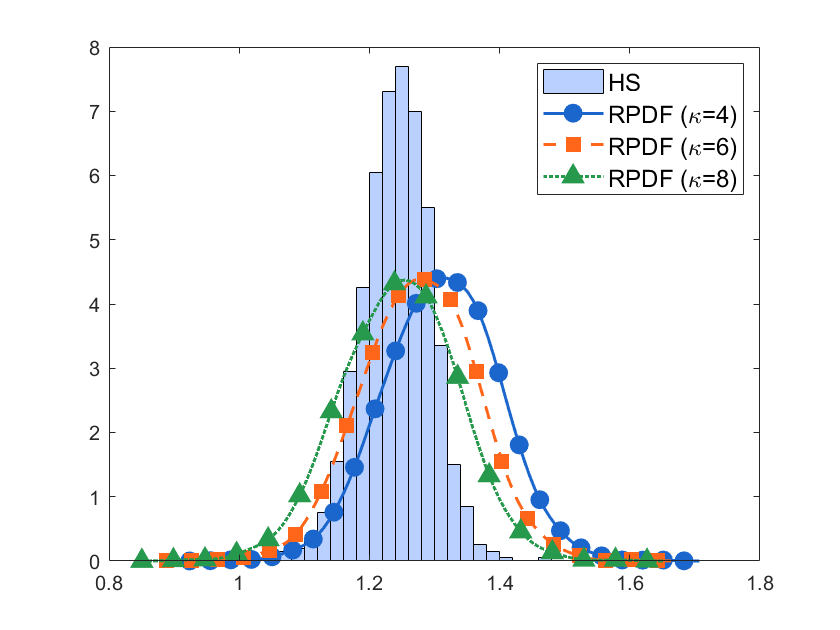}
\label{Fig4.2-1.3}
}
\subfigure[\,$x=2\pi$\,]{
\includegraphics[width=0.47\textwidth]{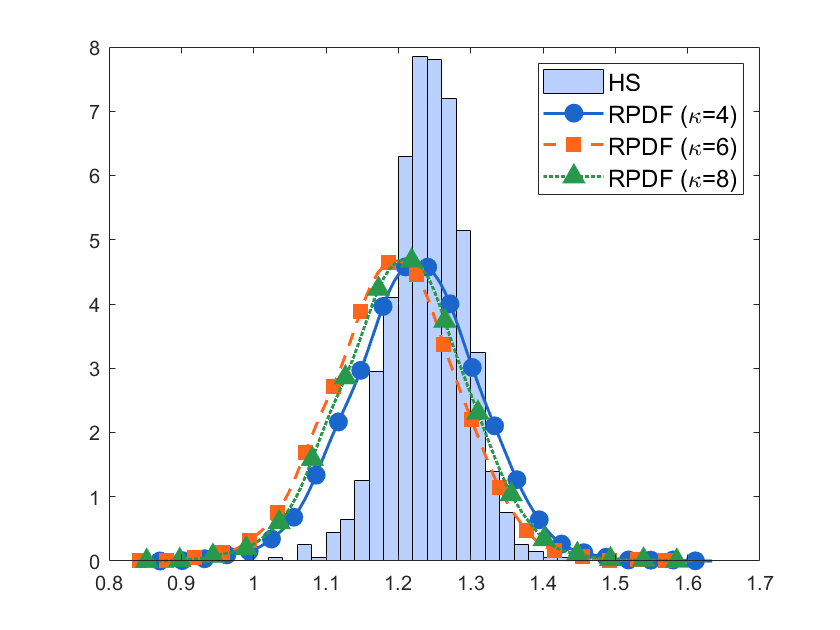}
\label{Fig4.2-1.4}
}
  \caption{Histograms of samples and results of reconstructed probability density function at four different positions for Example \ref{Example5}. }\label{Fig5.2}
\end{center}
\end{figure}

The exact covariance matrix is shown in Fig. \ref{Fig5.3}. We take the noise level $\tau=0.1\%$.
In order to show the effect of $\kappa_Q$ more clearly, we choose $\boldsymbol{\kappa}= [0.5, 1, 2, 3 ,4]^{\top}, [0.5, 1, 2, 3 ,4,5,6]^{\top}, [0.5, 1, 2, 3 ,4,5,6,7,8]^{\top}$.
The reconstruction results and errors are presented in Fig. \ref{Fig5.1} and Table \ref{tab5.1}, respectively.
We can find that the reconstruction is more satisfied as $\kappa_Q$ increases.
However, it should be noted that the relative error of the covariance matrix is a little large.
Besides, the Gibbs phenomenon occurs around the discontinuities of the surface, which drags down the overall reconstruction results.
In addition, by observing the reconstruction results of the probability density function in Fig. \ref{Fig5.2}, we can find an interesting phenomenon: Original generated sample distribution is right skewed ($S=0.3$), but reconstructed results are slightly left skewed.
Such a phenomenon means that the reconstructed results are somewhat biased.
It may be due to the slow decay of the Fourier coefficients of the binary grating, resulting in less accurate reconstruction results for each sample.
Inaccuracies in the reconstruction of individual samples can affect the estimation of the stochastic process's overall distribution, which in turn leads to some degree of deviation from the true distribution.

\begin{table}
\renewcommand{\arraystretch}{1.65}
\setlength{\tabcolsep}{15pt} 
\centering
\caption{The reconstruction errors of Example \ref{Example5} for $S=0.3$ and $K=4.0$.}\label{tab5.1}
\begin{tabular}{lll}
\Xhline{1.2pt}
$\kappa_Q$   &$Err_{\mathrm{mean}}$ & $Err_{\mathrm{cov}}$ \\ \hline
$4$ & $4.31\times 10^{-1}$ & 57.31\%      \\
$6$ & $2.43\times 10^{-1}$ & 55.58\%     \\
$8$ & $1.85\times 10^{-1}$ & 54.85\%   \\
\Xhline{1.2pt}
\end{tabular}
\end{table}

\section{Concluding remarks}
In this paper, we present an efficient numerical method for solving the inverse scattering problem of acoustic-elastic interaction with random periodic structures.
Our method reconstructs key statistical parameters of random structures using the acoustic scattered field data.
The mathematical model of the acoustic-elastic interaction are much more complicated than the problems with single acoustic or elastic wave field due to the coupling of different physical fields and the complex transmission conditions.
Since the reconstruction formula is explicit, the inversion process avoids repetitive computation of the direct problem, which significantly reduces the computational costs.
Thus, our algorithm is highly effective for acoustic-elastic interaction model.
Numerical examples demonstrate the effectiveness and reliability of the algorithm for stationary Gaussian and non-Gaussian stochastic processes.
In future work, we will extend our algorithm to other structures, such as obstacles \cite{NLvW2023MMAS} and cavities \cite{LLvW2024IPI}.

\flushleft{\textbf{Acknowledgements}}
This work of J.L. was partially supported by
the National Natural Science Foundation of China grant 12271209.

\section*{Declarations}
\textbf{Conﬂict of interest} The authors declare that they have no conﬂict of interest.
\vskip 5pt
\noindent
\textbf{Data availability} The data that support the findings of this study are available from the corresponding author upon reasonable request.



%
%



\end{document}